\documentclass{article}
\newcommand{\SAGAsection}{\textsc{Saga}}

\usepackage{amsmath}
\usepackage{amsthm}
\usepackage{amssymb}
\usepackage{empheq}
\usepackage{xcolor, color, colortbl}
\usepackage{mdframed}
\usepackage{wrapfig}
\usepackage{dsfont}
\usepackage{nccmath}
\usepackage{graphicx} %
\usepackage{subfigure}
\usepackage{natbib}

\usepackage{algorithm}
\usepackage{algorithmic}
\usepackage{stmaryrd}
\usepackage{nicefrac}
\usepackage{booktabs}
\usepackage{multirow}
\usepackage{rotating}
\usepackage{float}

\definecolor{mydarkblue}{rgb}{0,0.08,0.45}
\usepackage[colorlinks=true,
    linkcolor=mydarkblue,
    citecolor=mydarkblue,
    filecolor=mydarkblue,
    urlcolor=mydarkblue,
    pdfview=FitH,
    breaklinks]{hyperref}
\usepackage[final]{nips_2017}

\def\RR{{\mathbb R}}
\def\EE{{\mathbb E}}
\newcommand{\Econd}{\mathbf{E}}
\def\prox{{\mathbf{prox}}}
\def\eprox{{\emph{\text{prox}}}}

\newcommand*\mybluebox[1]{\colorbox{myblue}{\hspace{1em}#1\hspace{1em}}}

\newcommand{\PASAGA}{\textsc{ProxAsaga}}
\newcommand{\SAGA}{\textsc{Saga}}
\newcommand{\ASAGA}{\textsc{Asaga}}
\newcommand{\SAG}{\textsc{Sag}}

\newcommand{\SVRG}{\textsc{Svrg}}
\newcommand{\ProxSVRG}{Prox\textsc{Svrg}}
\newcommand{\SDCA}{\textsc{Sdca}}
\newcommand{\SGD}{\textsc{Sgd}}
\newcommand{\HOGWILD}{\textsc{Hogwild}}
\newcommand{\AsySPCD}{\textsc{AsySpcd}}

\DeclareMathOperator*{\argmin}{arg\,min}
\def\xx{{\boldsymbol x}}

\def\yy{{\boldsymbol y}}
\def\vv{{\boldsymbol v}}
\def\uu{{\boldsymbol u}}

\def\zz{{\boldsymbol z}}
\def\DD{{\boldsymbol D}}
\def\balpha{{\boldsymbol \alpha}}

\newcommand{\tablefont}[1] {{\fontsize{8}{10}\sffamily{#1}}}

\definecolor{myblue}{HTML}{D2E4FC}
\definecolor{Gray}{gray}{0.92}

\newmdtheoremenv{theo}{Theorem}
\newtheorem{theorem}{Theorem}
\newtheorem{lemma}{Lemma}
\newtheorem{definition}{Definition}
\newtheorem{corollary}{Corollary}
\newtheorem{innercustomtheorem}{Theorem}
\newenvironment{customtheorem}[1]
  {\renewcommand\theinnercustomtheorem{#1}\innercustomtheorem}
  {\endinnercustomtheorem}

\newtheorem{innercustomcorollary}{Corollary}
\newenvironment{customcorollary}[1]
  {\renewcommand\theinnercustomcorollary{#1}\innercustomcorollary}
  {\endinnercustomcorollary}

\graphicspath{{./figures/}}
\renewcommand{\llbracket}{[}
\renewcommand{\rrbracket}{]}

\title{Breaking the Nonsmooth Barrier: A Scalable Parallel Method for Composite Optimization}

\author{
    Fabian Pedregosa \\
    INRIA/ENS\thanks{DI \'{E}cole normale sup\'{e}rieure, CNRS, PSL Research University} \\ Paris, France\\
    \And
    R\'emi Leblond\\
    INRIA/ENS\footnotemark[1] \\
    Paris, France\\
    \And
    Simon Lacoste-Julien\\
    MILA and DIRO \\
    Universit\'{e} de Montr\'{e}al, Canada
}

\begin{document}
\maketitle

\begin{abstract}
  Due to their simplicity and excellent performance, parallel asynchronous variants of stochastic gradient descent have become popular methods to solve a wide range of large-scale optimization problems on multi-core architectures. Yet, despite their practical success, support for nonsmooth objectives is still lacking, making them unsuitable for many problems of interest in machine learning, such as the Lasso, group Lasso or empirical risk minimization with convex constraints.
  In this work, we propose and analyze \PASAGA, a fully asynchronous sparse method inspired by \SAGA, a variance reduced incremental gradient algorithm. The proposed method is easy to implement and significantly outperforms the state of the art on several nonsmooth, large-scale problems. We prove that our method achieves a theoretical linear speedup with respect to the sequential version under assumptions on the sparsity of gradients and block-separability of the proximal term. Empirical benchmarks on a multi-core architecture illustrate practical speedups of up to 12x on a 20-core machine.
\end{abstract}

\section{Introduction}
The widespread availability of multi-core computers motivates the development of parallel methods adapted for these architectures.
One of the most popular approaches is \HOGWILD~\citep{hogwild2011}, an asynchronous variant of stochastic gradient descent (\SGD).
In this algorithm, multiple threads run the update rule of \SGD\ asynchronously in parallel.
As \SGD, it only requires visiting a small batch of random examples per iteration, which makes it ideally suited for large scale machine learning problems.
Due to its simplicity and excellent performance, this parallelization approach has recently been extended to other variants of \SGD\ with better convergence properties, such as \SVRG~\citep{johnson2013accelerating} and \SAGA~\citep{defazio2014saga}.

Despite their practical success, existing parallel asynchronous variants of \SGD\ are limited to smooth objectives, making them inapplicable to many problems in machine learning and signal processing.
In this work, we develop a sparse variant of the \SAGA\ algorithm and consider its parallel asynchronous variants for general \emph{composite} optimization problems of the form:
\begin{empheq}[box=\mybluebox]{equation}\label{eq:opt}\tag{OPT}
\argmin^{\vphantom{i}}_{\xx \in \RR^p} f(\xx) \,+~ h(\xx) \quad,\quad \text{ with } f(\xx) := \textstyle\frac{1}{n}\sum_{i=1}^n f_i(\xx) \quad,
\end{empheq}
where each $f_i$ is convex with $L$-Lipschitz gradient, the average function $f$ is $\mu$-strongly convex and $h$ is convex but potentially nonsmooth.
We further assume that $h$ is ``simple'' in the sense that we have access to its proximal operator, and that it is
block-separable, that is, it can be decomposed block coordinate-wise as $h(\xx) = \textstyle{\sum_{B \in \mathcal{B}}} h_B(\llbracket\xx\rrbracket_B)$, where $\mathcal{B}$ is a partition of the coefficients into subsets which will call \emph{blocks} and $h_B$ only depends on coordinates in block $B$.  Note that there is no loss of generality in this last assumption as a unique block covering all coordinates is a valid partition, though in this case, our sparse variant of the \SAGA\ algorithm reduces to the original \SAGA\ algorithm and no gain from sparsity is obtained.

This template models a broad range of problems arising in machine learning and signal processing: the finite-sum structure of $f$ includes the least squares or logistic loss functions;
the proximal term $h$ includes penalties such as the $\ell_1$ or group lasso penalty.
Furthermore, this term can be extended-valued, thus allowing for convex constraints through the indicator function.

\paragraph{Contributions.} This work presents two main contributions. First, in \S\ref{scs:sparse_prox_saga} we describe Sparse Proximal \SAGA, a novel variant of the \SAGA\ algorithm which features a reduced cost per iteration in the presence of sparse gradients and a block-separable penalty.
Like other variance reduced methods, it enjoys a linear convergence rate under strong convexity.
Second, in \S\ref{sec:pasaga} we present \PASAGA, a lock-free asynchronous parallel version of the aforementioned algorithm that does not require consistent reads.
Our main results states that \PASAGA\ obtains (under assumptions) a theoretical linear speedup with respect to its sequential version.
Empirical benchmarks reported in \S\ref{scs:experiments} show that this method dramatically outperforms state-of-the-art alternatives on large sparse datasets, while the empirical speedup analysis illustrates the practical gains as well as its limitations.

\subsection{Related work}

\paragraph{Asynchronous coordinate-descent.} For composite objective functions of the form~\eqref{eq:opt}, most of the existing literature on asynchronous optimization has focused on variants of coordinate descent.
\citet{liu2015asynchronous2} proposed an asynchronous variant of (proximal) coordinate descent and proved a near-linear speedup in the number of cores used, given a suitable step size.
This approach has been recently extended to general block-coordinate schemes by \citet{peng2016arock}, to greedy coordinate-descent schemes by~\citet{you2016asynchronous} and to non-convex problems by~\citet{apalm_davis}. However, as illustrated by our experiments, in the large sample regime coordinate descent compares poorly against incremental gradient methods like \SAGA.

\paragraph{Variance reduced incremental gradient and their asynchronous variants.} Initially proposed in the context of smooth optimization by \citet{roux2012stochastic}, variance reduced incremental gradient methods have since been extended to minimize composite problems of the form~\eqref{eq:opt} (see table below). Smooth variants of these methods have also recently been extended to the asynchronous setting, where multiple threads run the update rule asynchronously and in parallel.
Interestingly,
none of these methods achieve both simultaneously, i.e. asynchronous optimization of composite problems.
Since variance reduced incremental gradient methods have shown state of the art performance in both settings, this generalization is of key practical interest.
\begin{table}[H]
\centering\footnotesize
\hskip-0.5cm\begin{tabular}{c c | l l|}
\cline{2-4}
\multicolumn{1}{c|}{} & {\tablefont{Objective}} & {\hskip 1.5cm \tablefont{Sequential Algorithm}} & {\hskip 1.5cm \tablefont{Asynchronous Algorithm}}\\
\cline{2-4}
\multicolumn{1}{c|}{\multirow{0}{*}{}}
& \cellcolor{Gray} &\cellcolor{Gray}{\SVRG~\citep{johnson2013accelerating}}&
 \cellcolor{Gray}{\SVRG~\citep{reddi2015variance}}\\
\multicolumn{1}{c|}{} &\multirow{-1}{*}{\tablefont{Smooth}}\cellcolor{Gray}&\SDCA~\citep{shalev2013stochastic}
 \cellcolor{Gray} &
\cellcolor{Gray}\textsc{Passcode}~\citep[\SDCA\ variant]{hsieh2015passcode}
\\
\multicolumn{1}{c|}{} &\cellcolor{Gray} &\SAGA~\citep{defazio2014saga}
 \cellcolor{Gray} &\cellcolor{Gray}\ASAGA~\citep[\SAGA\ variant]{leblond2016Asaga}\\
\multicolumn{1}{c|}{} &\multirow{0}{*}{}  &\textsc{ProxSdca}~\citep{shalev2012proximal} &  \\
\multicolumn{1}{c|}{}& \multirow{-1}{*}{\tablefont{Composite}} & 
\SAGA~\citep{defazio2014saga}  & {\qquad\qquad {This work: \PASAGA}}\\
\multicolumn{1}{c|}{}& \multirow{-2}{*}{} &
\ProxSVRG~\citep{xiao2014proximal} &  \\
\cline{2-4}\\
\end{tabular}
\end{table}

\paragraph{On the difficulty of a composite extension.}
Two key issues explain the paucity in the development of asynchronous incremental gradient methods for composite optimization.
The first issue is related to the design of such algorithms.
Asynchronous variants of \SGD\ are most competitive when the updates are sparse and have a small overlap, that is, when each update modifies a small and different subset of the coefficients. This is typically achieved by updating only coefficients for which the partial gradient at a given iteration is nonzero,\footnote{Although some regularizers are sparsity inducing, large scale datasets are often extremely sparse and leveraging this property is crucial for the efficiency of the method.} but existing schemes such as the lagged updates technique~\citep{schmidt2016minimizing} are not applicable in the asynchronous setting.
The second difficulty is related to the analysis of such algorithms.
All convergence proofs crucially use the Lipschitz condition on the gradient to bound the noise terms derived from asynchrony.
However, in the composite case, the gradient mapping term~\citep{beck2009gradient}, which replaces the gradient in proximal-gradient methods, does not have a bounded Lipschitz constant.
Hence, the traditional proof technique breaks down in this scenario.

\paragraph{Other approaches.}
Recently, \citet{meng2017aaai, gu2016asynchronous} independently proposed a doubly stochastic method to solve the problem at hand. Following \citet{meng2017aaai} we refer to it as Async-\textsc{ProxSvrcd}.
This method performs coordinate descent-like updates in which the true gradient is replaced by its \SVRG\ approximation.
It hence features a doubly-stochastic loop: at each iteration we select a random coordinate \emph{and} a random sample.
Because the selected coordinate block is uncorrelated with the chosen sample, the algorithm can be orders of magnitude slower than \SAGA\ in the presence of sparse gradients. \ref{apx:experiments} contains a comparison of these methods.

\subsection{Definitions and notations}
By convention, we denote vectors and vector-valued functions in lowercase boldface (e.g. $\xx$) and matrices in uppercase boldface  (e.g. $\boldsymbol D$).
The { proximal operator } of a convex lower semicontinuous function $h$ is defined as
$  \prox_{ h}(\xx) := \argmin_{\zz \in \RR^p} \{ h(\zz) + \frac{1}{2}\|\xx - \zz\|^2\}$.
A function $f$ is said to be {$L$-smooth} if it is differentiable and its gradient is $L$-Lipschitz continuous.
A function $f$ is said to be {$\mu$-strongly convex} if $f - \frac{\mu}{2}\|\cdot\|^2$ is convex.
We use the notation $\kappa := L/\mu$ to denote the condition number for an $L$-smooth and $\mu$-strongly convex function.\footnote{Since we have assumed that each individual $f_i$ is $L$-smooth, $f$ itself is $L$-smooth -- but it could have a smaller smoothness constant. Our rates are in terms of this bigger $L/\mu$, as is standard in the \SAGA\ literature.}

$\boldsymbol{I}_{p}$ denotes the $p$-dimensional identity matrix, $\mathds{1}\{\text{cond}\}$ the characteristic function, which is $1$ if cond evaluates to true and $0$ otherwise.
The average of a vector or matrix is denoted $\overline{\boldsymbol\alpha}$ $:= \frac{1}{n}\sum_{i=1}^n {\boldsymbol\alpha}_i$.
We use $\| \cdot \|$ for the Euclidean norm.
For a positive semi-definite matrix~$\DD$, we define its associated distance as $\|\xx\|^2_{\boldsymbol{D}} := \langle \xx, \boldsymbol{D} \xx\rangle$. We denote by $[\,\xx\,]_b$ the $b$-th coordinate in $\xx$. This notation is overloaded so that for a collection of blocks $T = \{B_1, B_2, \ldots\}$, $\llbracket \xx \rrbracket_T$ denotes the vector $\xx$ restricted to the coordinates in the blocks of $T$. For convenience, when $T$ consists of a single block $B$ we use $\llbracket \xx \rrbracket_B$  as a shortcut of $\llbracket \xx \rrbracket_{\{B\}}$.
Finally, we distinguish $\EE$, the full expectation taken with respect to all the randomness in the system, from $\Econd$, the conditional expectation of a random $i_t$ (the random feature sampled at each iteration by \SGD-like algorithms) conditioned on all the ``past'', which the context will clarify.

\section{Sparse Proximal SAGA }\label{scs:sparse_prox_saga}

\paragraph{Original \SAGAsection\ algorithm.} The original \SAGA\ algorithm~\citep{defazio2014saga} maintains two moving quantities: the current iterate~$\xx$ and a table (memory) of historical gradients $({\boldsymbol\alpha}_i)_{i=1}^n$.
At every iteration, it samples an index $i \in \{1,\ldots, n\}$ uniformly at random, and computes the next iterate $(\xx^+, \balpha^+)$ according to the following recursion:
\begin{equation}\label{eq:original_saga}
  \uu_i = \nabla f_i(\xx) - {\boldsymbol\alpha}_i + \overline{\boldsymbol\alpha} \,;\quad\xx^+ = \prox_{\gamma h}\big(\xx - \gamma \uu_i\big);\quad {\boldsymbol\alpha}_i^+ = \nabla f_i(\xx)~.
\end{equation}
On each iteration, this update rule requires to visit all coefficients even if the partial gradients $\nabla f_i$ are sparse.
Sparse partial gradients arise in a variety of practical scenarios: for example, in generalized linear models the partial gradients inherit the sparsity pattern of the dataset.
Given that large-scale datasets are often sparse,\footnote{For example, in the \texttt{LibSVM} datasets suite, 8 out of the 11 datasets (as of May 2017) with more than a million samples have a density between $10^{-4}$ and $10^{-6}$.} leveraging this sparsity is crucial for the success of the optimizer.

\paragraph{Sparse Proximal \SAGAsection\ algorithm.} We will now describe an algorithm that leverages sparsity in the partial gradients by only updating those blocks that intersect with the support of the partial gradients. 
Since in this update scheme some blocks might appear more frequently than others, we will need to counterbalance this undersirable effect with a well-chosen block-wise reweighting of the average gradient and the proximal term. 

In order to make precise this block-wise reweighting, we define the following quantities.
We denote by $T_i$ the \emph{extended support} of $\nabla f_i$, which is the set of blocks that intersect the support of $\nabla f_i$, formally defined as $T_i := \{B: \text{supp}(\nabla f_i) \cap B \neq \varnothing, \,B\in\mathcal{B} \}$.
For totally separable penalties such as the $\ell_1$ norm, the blocks are individual coordinates and so the extended support covers the same coordinates as the support.
Let ${d}_B := n / n_B$, where $n_B:=\sum_i\mathds{1}\{B \in T_i\}$ is the number of times that $B \in T_i$. For simplicity we assume $n_B > 0$, as otherwise the problem can be reformulated without block $B$.
The update rule in \eqref{eq:original_saga} requires computing the proximal operator of $h$, which involves a full pass on the coordinates.
In our proposed algorithm, we replace~$h$ in~\eqref{eq:original_saga} with the function $\varphi_i(\xx) := \textstyle\sum_{B \in T_{i}} d_B h_B(\xx)$, whose form is justified by the following
 three properties. First, this function is zero outside $T_i$, allowing for sparse updates. Second, because of the block-wise reweighting~$d_B$, the function $\varphi_i$ is an unbiased estimator of $h$ (i.e.,  $\Econd\, \varphi_i = h$), property which will be crucial to prove the convergence of the method.
Third, $\varphi_i$ inherits the block-wise structure of $h$ and its proximal operator can be computed from that of $h$ as 
$
\llbracket\prox_{\gamma \varphi_i}(\xx)\rrbracket_{B} = \llbracket\prox_{(d_B \gamma) h_B}(\xx)\rrbracket_B$ if $B \in T_i$ and $
\llbracket\prox_{\gamma \varphi_i}(\xx)\rrbracket_{B} = \llbracket\xx\rrbracket_B$ otherwise.
Following~\citet{leblond2016Asaga}, we will also replace the dense gradient estimate $\uu_i$ by the sparse estimate $\vv_i := \nabla f_i(\xx) - \balpha_i + \DD_i \overline{\balpha}$, where $\DD_i$ is the diagonal matrix defined block-wise as $\llbracket\DD_i\rrbracket_{B, B} = d_B \mathds{1}\{B \in T_i\}\boldsymbol{I}_{|B|}$. 
It is easy to verify that the vector $\DD_i \overline{\balpha}$ is a weighted projection onto the support of $T_i$ and
$\Econd\,\DD_i \overline{\balpha} =\overline{\balpha}$,
making $\vv_i$ an unbiased estimate of the gradient.

We now have all necessary elements to describe the Sparse Proximal \SAGA\ algorithm. As the original \SAGA\ algorithm, it maintains two moving quantities: the current iterate~$\xx \in \RR^p$ and a table of historical gradients $({\boldsymbol\alpha}_i)_{i=1}^n,\, \balpha_i \in \RR^p$.
At each iteration, the algorithm samples an index $i \in \{1, \ldots, n\}$ and computes the next iterate $(\xx^+, \balpha^+)$ as:
\begin{empheq}[box=\mybluebox]{equation}\label{eq:SPS}\tag{SPS}
    \vphantom{\sum_i^n}\vv_i = \nabla f_i(\xx) - \balpha_i + \DD_i \overline{\balpha}\,; ~\xx^+ = \prox_{\gamma \varphi_i}\big(\xx - \gamma \vv_i \big)\,;~\balpha_i^+ = \nabla f_i(\xx)~,
\end{empheq}
where in a practical implementation the vector $\overline{\boldsymbol\alpha}$ is updated incrementally at each iteration.

The above algorithm is sparse in the sense that it only requires to visit and update blocks in the extended support: if $B \notin T_i$, by the sparsity of $\vv_i$ and $\prox_{\varphi_i}$, we have $\llbracket\xx^+\rrbracket_B = \llbracket\xx\rrbracket_B$.
Hence, when the extended support $T_i$ is sparse, this algorithm can be orders of magnitude faster than the naive \SAGA\ algorithm. The extended support is sparse for example when the partial gradients are sparse and the penalty is separable, as is the case of the $\ell_1$ norm or the indicator function over a hypercube, or when the the penalty is block-separable in a way such that only a small subset of the blocks overlap with the support of the partial gradients.
Initialization of variables and a reduced storage scheme for the memory are discussed in the implementation details section of \ref{apx:implementation_details}.

{\bfseries Relationship with existing methods}.
This algorithm can be seen as a generalization of both the Standard \SAGA\ algorithm and the Sparse \SAGA\ algorithm of \citet{leblond2016Asaga}.
When the proximal term is not block-separable, then $d_B=1$ (for a unique block $B$) and the algorithm defaults to the Standard (dense) \SAGA\ algorithm.
In the smooth case (i.e., $h=0$), the algorithm defaults to the Sparse \SAGA\ method.
Hence we note that the sparse gradient estimate $\vv_i$ in our algorithm is the same as the one proposed in~\citet{leblond2016Asaga}.
However, we emphasize that a straightforward combination of this sparse update rule with the proximal update from the Standard \SAGA\ algorithm results in a nonconvergent algorithm: the block-wise reweighting of $h$ is a surprisingly simple but crucial change.
We now give the convergence guarantees for this algorithm.

\begin{theorem}\label{theorem:rates_sparse_saga}
  \label{th1} Let $\gamma = \frac{a}{5L}$ for any $a\leq 1$ and $f$ be $\mu$-strongly convex ($\mu > 0$). Then Sparse Proximal \SAGA\ converges geometrically in expectation with a rate factor of at least $\rho = \frac{1}{5} \min\{\frac{1}{n}, a\frac{1}{\kappa}\}$. That is, for $\xx_t$ obtained after $t$ updates, we have the following bound:
  $$
  \EE \|\xx_t - \xx^*\|^2 \leq (1-\rho)^t C_0\,,~\text{ with }C_0 := \|\xx_0 - \xx^*\|^2 + \textstyle\frac{1}{5 L^2} \textstyle\sum_{i=1}^n\|{\boldsymbol\alpha}_i^0 - \nabla f_i(\xx^*)\|^2  \, \quad.
  $$
\end{theorem}
{\bfseries Remark}. For the step size $\gamma=\nicefrac{1}{5L}$, the convergence rate is $(1 - \nicefrac{1}{5}\min\{\nicefrac{1}{n}, \nicefrac{1}{\kappa}\})$.
We can thus identify two regimes: the ``big data'' regime, $n \geq \kappa$, in which the rate factor is bounded by $ \nicefrac{1}{5n}$, and the ``ill-conditioned'' regime, $\kappa \geq n$, in which the rate factor is bounded by $ \nicefrac{1}{5\kappa}$.
This rate roughly matches the rate obtained by \citet{defazio2014saga}.
While the step size bound of $\nicefrac{1}{5L}$ is slightly smaller than the $\nicefrac{1}{3L}$ one obtained in that work, this can be explained by their stronger assumptions: each $f_i$ is strongly convex whereas they are strongly convex only on average in this work. All proofs for this section can be found in \ref{apx:sparse_saga}.

\section{Asynchronous Sparse Proximal SAGA}\label{sec:pasaga}

We introduce \PASAGA\ -- the asynchronous parallel variant of Sparse Proximal \SAGA. In this algorithm, multiple cores update a central parameter vector using the Sparse Proximal \SAGA\ introduced in the previous section, and updates are performed asynchronously. The algorithm parameters are read and written without vector locks, i.e., the vector content of the shared memory can potentially change while a core is reading or writing to main memory coordinate by coordinate. These operations are typically called \emph{inconsistent} (at the vector level).

The full algorithm is described in
Algorithm~\ref{alg:theoretical} for its theoretical version (on which our analysis is built) and in Algorithm~\ref{alg:sagasync} for its practical implementation. The practical implementation differs from the analyzed agorithm in three points.
First, in the implemented algorithm, index $i$ is sampled before reading the coefficients to minimize memory access since only the extended support needs to be read. Second, since our implementation targets generalized linear models, the memory $\boldsymbol\alpha_i$ can be compressed into a single scalar in L\ref{line:alphaiupdate} (see \ref{apx:implementation_details}). Third, $\overline\balpha$ is stored in memory and updated incrementally instead of recomputed at each iteration. %

The rest of the section is structured as follows:
we start by describing our framework of analysis;
we then derive essential properties of \PASAGA\ along with a classical delay assumption.
Finally, we state our main convergence and speedup result.

\begin{figure*}[ttt!]
 \begin{minipage}[t]{0.5\textwidth}
   \begin{algorithm}[H]
     \caption{\PASAGA\ (analyzed)}
     \label{alg:theoretical}
     \label{theoreticalgo}
     \begin{algorithmic}[1]
    \STATE Initialize shared variables $\xx$ and $({\boldsymbol\alpha}_i)_{i=1}^n$
    \LOOP
    \STATE $\hat \xx = $ inconsistent read of $\xx$
 	  \STATE $\hat {\boldsymbol\alpha} = $ inconsistent read of ${\boldsymbol\alpha}$
    	  \STATE \emph{Sample} $i$  uniformly in $\{1,...,n\}$
        \STATE $S_i :=$ support of $\nabla f_i$
        \STATE $T_i :=$ extended support of $\nabla f_i$ in $\mathcal{B}$
        \STATE $\llbracket\,\overline{{\boldsymbol\alpha}}\,\rrbracket_{T_i} = \nicefrac{1}{n} \sum_{j=1}^n \llbracket\,\hat{{\boldsymbol\alpha}}_j\,\rrbracket_{T_i}$
        \STATE $\llbracket\,\delta\balpha\,\rrbracket_{S_i} = \llbracket\nabla f_i(\hat \xx)\rrbracket_{S_i} - \llbracket\hat{{\boldsymbol\alpha}}_i\rrbracket_{S_i}$
        \STATE $\llbracket\,\hat \vv\,\rrbracket_{T_i} = \llbracket\,\delta\balpha\,\rrbracket_{T_i} + \llbracket\DD_i \overline {\boldsymbol\alpha}\,\rrbracket_{T_i}$
        \STATE $\llbracket\,\delta\xx\,\rrbracket_{T_i} = \llbracket\prox_{\gamma \varphi_i}(\hat{\xx} - \gamma \hat \vv)\rrbracket_{T_i} - \llbracket\hat\xx\rrbracket_{T_i}$
        \FOR{$B$ \,{\bfseries in} $T_i$}
        \FOR{$b \in B$}
        \STATE $[\,\xx\,]_b \gets [\,\xx\,]_b + [\,\delta\xx\,]_b$\hfill$\triangleright$ atomic \label{alg1:atomicLine}
        \IF{$b \in S_i$}
        \STATE $[{\boldsymbol\alpha}_i]_b \gets [\nabla f_i(\hat{\xx})]_b$
        \ENDIF
        \ENDFOR
        \ENDFOR
        \STATE // {\footnotesize(`$\gets$' denotes shared memory update.)}
    \ENDLOOP
   \end{algorithmic}
    \end{algorithm}
 \end{minipage}
 \hfill
 \begin{minipage}[t]{0.5\textwidth}
    \begin{algorithm}[H]
      \caption{\PASAGA\ (implemented)}
      \label{alg:sagasync}
      \begin{algorithmic}[1]
     \STATE Initialize shared variables $\xx$, $({\boldsymbol\alpha}_i)_{i=1}^n$, $\overline{\balpha}$
     \LOOP
     \STATE \emph{Sample} $i$ uniformly in $\{1,...,n\}$
     \STATE $S_i :=$ support of $\nabla f_i$
     \STATE $T_i :=$ extended support of $\nabla f_i$ in $\mathcal{B}$
     \STATE $\llbracket\,\hat \xx\,\rrbracket_{T_i} = $ inconsistent read of $\xx$ on $T_i$
     \STATE $\hat {\boldsymbol\alpha}_i = $ inconsistent read of ${\boldsymbol\alpha}_i$
     \STATE $\llbracket\,\overline{\boldsymbol\alpha}\,\rrbracket_{T_i} = $ inconsistent read of $\overline{\boldsymbol\alpha}$ on $T_i$
     \STATE $\llbracket\,\delta\balpha\,\rrbracket_{S_i} = \llbracket\nabla f_i(\hat \xx)\rrbracket_{S_i} - \llbracket\hat{{\boldsymbol\alpha}}_i\rrbracket_{S_i}$
     \STATE $\llbracket\,\hat \vv\,\rrbracket_{T_i} = \llbracket\delta \balpha\,\rrbracket_{T_i} + \llbracket\,\DD_i \overline {\boldsymbol\alpha}\,\rrbracket_{T_i}$
     \STATE $\llbracket\,\delta\xx\,\rrbracket_{T_i} = \llbracket\prox_{\gamma \varphi_i}(\hat{\xx} - \gamma \hat \vv)\rrbracket_{T_i} - \llbracket\hat\xx\rrbracket_{T_i}$
     \FOR{$B$ \,{\bfseries in} $T_i$}
     \FOR{$b$ \,{\bfseries in} $B$}
     \STATE $[\,\xx\,]_b \gets [\,\xx\,]_b + [\,\delta\xx\,]_b$\hfill$\triangleright$ atomic
     \IF{$b \in S_i$}
      \STATE $[\,\overline {\boldsymbol\alpha}\,]_b \gets [\overline {\boldsymbol\alpha}]_b + \nicefrac{1}{n}[\delta {\boldsymbol\alpha}]_b$\hfill$\triangleright$ atomic
     \ENDIF
     \ENDFOR
     \ENDFOR
     \STATE ${\boldsymbol\alpha}_i \gets \nabla f_i(\hat\xx)$  \quad  (scalar update)\hfill$\triangleright$ atomic\label{line:alphaiupdate}
     \ENDLOOP
      \end{algorithmic}
    \end{algorithm}
 \end{minipage}
\end{figure*}

\subsection{Analysis framework}
As in most of the recent asynchronous optimization literature, we build on the hardware model introduced by~\citet{hogwild2011}, with multiple cores reading and writing to a shared memory parameter vector.
These operations are asynchronous (lock-free) and \textit{inconsistent}:\footnote{This is an extension of the framework of \citet{hogwild2011}, where consistent updates were assumed.} $\hat \xx_t$, the local copy of the parameters of a given core, does not necessarily correspond to a consistent iterate in memory.

\paragraph{``Perturbed'' iterates.}
To handle this additional difficulty, contrary to most contributions in this field, we choose the ``perturbed iterate framework'' proposed by~\citet{mania2015perturbed} and refined by~\citet{leblond2016Asaga}.
This framework can  analyze variants of \SGD\ which obey the update rule:
\begin{equation*}
\xx_{t+1} = \xx_t - \gamma \vv(\xx_t, i_t)\, ,~\text{ where $\vv$ verifies the unbiasedness condition $\Econd\, \vv(\xx, i_t) = \nabla f(\xx)$}
\end{equation*}
and the expectation is computed with respect to $i_t$.
In the asynchronous parallel setting, cores are reading inconsistent iterates from memory, which we denote $\hat{\xx}_t$.
As these inconsistent iterates are affected by various delays induced by asynchrony, they cannot easily be written as a function of their previous iterates.
To alleviate this issue,~\citet{mania2015perturbed} choose to introduce an additional quantity for the purpose of the analysis:
\begin{equation}
\xx_{t+1} := \xx_t - \gamma \vv(\hat \xx_t, i_t)\,, \quad\text{ the ``virtual iterate'' -- which is never actually computed}\, .
\end{equation}
Note that this equation is the \textit{definition} of this new quantity~$\xx_t$.
This virtual iterate is useful for the convergence analysis and makes for much easier proofs than in the related literature.

\paragraph{``After read'' labeling.}
How we choose to define the iteration counter $t$ to label an iterate $\xx_t$ matters in the analysis.
In this paper, we follow the ``after read'' labeling proposed in \citet{leblond2016Asaga}, in which
we update our iterate counter, $t$, as each core \emph{finishes reading} its copy of the parameters (in the specific case of \PASAGA, this includes both~$\hat \xx_t$ and~$\hat {\boldsymbol\alpha}^t$).
This means that $\hat \xx_t$ is the $(t+1)^{th}$ fully completed read.
One key advantage of this approach compared to the classical choice of~\citet{hogwild2011} -- where $t$ is increasing after each successful update -- is that it guarantees both that the $i_t$ are uniformly distributed and that~$i_t$ and~$\hat \xx_t$ are independent.
This property is not verified when using the ``after write'' labeling of~\citet{hogwild2011}, although it is still implicitly assumed in the papers using this approach, see \citet[Section 3.2]{leblond2016Asaga} for a discussion of issues related to the different labeling schemes.

\paragraph{Generalization to composite optimization.}
Although the perturbed iterate framework was designed for gradient-based updates, we can extend it to proximal methods by remarking that in the sequential setting, proximal stochastic gradient descent and its variants can be characterized by the following similar update rule:
\begin{equation}\label{eq:def_grad_mapping}
\xx_{t+1} = \xx_t - \gamma {\boldsymbol g}(\xx_t, \vv_{i_t}, {i_t}) \,,\quad \text{ with }\, {\boldsymbol g}(\xx, \vv, i) := \textstyle\frac{1}{\gamma}\big(\xx - \prox_{\gamma \varphi_i}(\xx - \gamma \vv)\big)\,,
\end{equation}
where as before $\vv$ verifies the unbiasedness condition $\Econd\, \vv = \nabla f(\xx)$.
The Proximal Sparse \SAGA\ iteration can be easily written within this template by using $\varphi_{i}$ and $\vv_i$ as defined in \S\ref{scs:sparse_prox_saga}. Using this definition of $\boldsymbol g$,
we can define \PASAGA\ virtual iterates as:
\begin{equation}\label{eq:definition}
\xx_{t+1} := \xx_t - \gamma \boldsymbol g(\hat{\xx}_t, \hat\vv^t_{i_t}, {i_t}) \,,\quad\text{ with } \hat{\vv}_{i_t}^t = \nabla f_{i_t}(\hat \xx_t) - \hat{{\boldsymbol\alpha}}^t_{i_t} + {\boldsymbol D}_{i_t} \overline{{\boldsymbol\alpha}}^t\quad,
\end{equation}
 where as in the sequential case, the memory terms are updated as
$\hat {\boldsymbol\alpha}^t_{i_t} = \nabla f_{i_t}(\hat \xx_t)$. Our theoretical analysis of \PASAGA\ will be based on this definition of the virtual iterate $\xx_{t+1}$.

\subsection{Properties and assumptions}
Now that we have introduced the ``after read'' labeling for proximal methods in Eq.~\eqref{eq:definition}, we can leverage the framework of~\citet[Section 3.3]{leblond2016Asaga} to derive essential properties for the analysis of \PASAGA.
We describe below three useful properties arising from the definition of Algorithm~\ref{alg:theoretical}, and then state a central (but standard) assumption that the delays induced by the asynchrony are uniformly bounded. 

{\bfseries Independence:} Due to the ``after read'' global ordering, $i_r$ is independent of $\hat{\xx}_t$ for all $r \geq t$. We enforce the independence for $r = t$ by having the cores read all the shared parameters before their iterations.

{\bfseries Unbiasedness:} The term $\hat \vv^t_{i_t}$ is an unbiased estimator of the gradient of $f$ at $\hat \xx_t$.
This property is a consequence of the independence between $i_t$ and $\hat{\xx}_t$.

{\bfseries Atomicity:} The shared parameter coordinate update of $[\xx]_b$ on Line~\ref{alg1:atomicLine} is atomic.
This means that there are no overwrites for a single coordinate even if several cores compete for the same resources. Most modern processors have support for atomic operations with minimal overhead.

{\bfseries Bounded overlap assumption.} We assume that there exists a uniform bound, $\tau$, on the maximum number of overlapping iterations.
This means that every coordinate update from iteration $t$ is successfully written to memory before iteration $t + \tau + 1$ starts.
Our result will give us conditions on $\tau$ to obtain linear speedups.

{\bfseries Bounding $\hat \xx_t - \xx_t$.}
The delay assumption of the previous paragraph allows to express the difference between real and virtual iterate using the gradient mapping  $\boldsymbol g_u := \boldsymbol g(\hat{\xx}_u, \hat\vv^u_{i_u}, {i_u})$ as:
\begin{equation}
\hat \xx_t - \xx_t = \gamma \textstyle\sum_{u=(t - \tau)_+}^{t-1}\boldsymbol G_{u}^t \boldsymbol g_u \,, \text{where $\boldsymbol G^t_{u}$ are $p\times p$ diagonal matrices with terms in $\{0, +1\}$.}
\end{equation}
$0$ represents instances where both $\hat \xx_u$ and $\xx_u$ have received the corresponding updates.
$+1$, on the contrary, represents instances where $\hat \xx_u$ has not yet received an update that is already in $\xx_u$ by definition. This bound will prove essential to our analysis.

\subsection{Analysis}\label{ssec:analysis}
In this section, we state our convergence and speedup results for \PASAGA.
The full details of the analysis can be found in \ref{apx:async}.
Following~\citet{hogwild2011}, we introduce a sparsity measure (generalized to the composite setting) that will appear in our results.
\begin{definition}
	Let $\Delta := \max_{B \in \mathcal{B}} |\{i: T_i \ni B\}| / n$.
	This is the normalized maximum number of times that a block appears in the extended support. For example, if a block is present in all $T_i$, then $\Delta=1$. If no two $T_i$ share the same block, then $\Delta = \nicefrac{1}{n}$. We always have $\nicefrac{1}{n} \leq \Delta \leq 1$.
\end{definition}

\begin{theorem}[Convergence guarantee of \PASAGA]\label{thm:convergence}
	Suppose $\tau \leq \frac{1}{10 \sqrt{\Delta}}$.
	For any step size $\gamma = \frac{a}{L}$ with $a \leq a^*(\tau) := \frac{1}{36} \min\{1, \frac{6 \kappa}{\tau}\}$, the inconsistent read iterates of Algorithm~\ref{alg:theoretical} converge in expectation at a geometric rate factor of at least: $\rho(a) = \frac{1}{5} \min \big\{\frac{1}{n},  a \frac{1}{\kappa}\big\},$
	i.e. $\EE \|\hat{\xx}_t-\xx^*\|^2 \leq (1-\rho)^t \,  \tilde C_0$, where $\tilde C_0$ is a constant independent of $t$ ($\approx \frac{n \kappa}{a}C_0$ with $C_0$ as defined in Theorem~\ref{th1}).
\end{theorem}
This last result is similar to the original \SAGA\ convergence result and our own Theorem~\ref{th1}, with both an extra condition on $\tau$ and on the maximum allowable step size. In the best sparsity case, $\Delta = \nicefrac{1}{n}$ and we get the condition $\tau \leq \nicefrac{\sqrt{n}}{10}$.
We now compare the geometric rate above to the one of Sparse Proximal \SAGA\ to derive the necessary conditions under which \PASAGA\ is linearly faster.

\begin{corollary}[Speedup]\label{thm:corollary}
	Suppose $\tau \leq \frac{1}{10 \sqrt{\Delta}}$.
	If $\kappa \geq n$, then using the step size $\gamma = \nicefrac{1}{36 L}$, \PASAGA\ converges geometrically with rate factor $\Omega(\frac{1}{\kappa})$.
	If $\kappa < n$, then using the step size $\gamma = \nicefrac{1}{36 n \mu}$, \PASAGA\ converges geometrically with rate factor $\Omega(\frac{1}{n})$.
	In both cases, the convergence rate is the same as Sparse Proximal \SAGA. Thus \PASAGA\ is linearly faster than its sequential counterpart up to a constant factor. Note that in both cases \emph{the step size does not depend on~$\tau$}.

	Furthermore, if $\tau \leq 6 \kappa$, we can use a universal step size of $\Theta(\nicefrac{1}{L})$ to get a similar rate for \PASAGA\ than Sparse Proximal \SAGA, thus making it adaptive to local strong convexity since the knowledge of $\kappa$ is not required.
\end{corollary}

These speedup regimes are comparable with the best ones obtained in the smooth case, including~\citet{hogwild2011,reddi2015variance}, even though unlike these papers, we support inconsistent reads and nonsmooth objective functions.
The one exception is~\citet{leblond2016Asaga}, where the authors prove that their algorithm, \ASAGA, can obtain a linear speedup even without sparsity in the well-conditioned regime.
In contrast, \PASAGA\ always requires some sparsity.
Whether this property for smooth objective functions could be extended to the composite case remains an open problem.

Relative to \AsySPCD, in the best case scenario (where the components of the gradient are uncorrelated, a somewhat unrealistic setting), \AsySPCD\ can get a near-linear speedup for~$\tau$ as big as~$\sqrt[4]{p}$.
Our result states that $\tau = \mathcal{O}(\nicefrac{1}{\sqrt{\Delta}})$ is necessary for a linear speedup.
This means in case $\Delta \leq \nicefrac{1}{\sqrt{p}}$ our bound is better than the one obtained for \AsySPCD.
Recalling that $\nicefrac{1}{n} \leq \Delta \leq 1$, it appears that \PASAGA\ is favored when $n$ is bigger than $\sqrt{p}$ whereas \AsySPCD\ may have a better bound otherwise, though this comparison should be taken with a grain of salt given the assumptions we had to make to arrive at comparable quantities.
An extended comparison with the related work can be found in \ref{apx:related_work}.

\section{Experiments} \label{scs:experiments}

In this section, we compare \PASAGA\ with related methods on different datasets.
Although \PASAGA\ can be applied more broadly, we focus on $\ell_1\!+\!\ell_2$-regularized logistic regression, a model of particular practical importance.
The objective function takes the form
\begin{equation}\label{eq:logistic_loss}
  \frac{1}{n} \sum_{i=1}^n \log\big(1 + \exp(- b_i \boldsymbol a_i^\intercal \xx)\big) + \textstyle\frac{\lambda_1}{2} \|\xx\|_2^2 + \lambda_2 \|\xx\|_1\quad,
\end{equation}
where $\boldsymbol a_i \in \mathbb{R}^p$ and $b_i \in \{-1,+1\}$ are the data samples. Following~\citet{defazio2014saga}, we set $\lambda_1 = 1/n$. The amount of $\ell_1$ regularization ($\lambda_2$) is selected to give an approximate $\nicefrac{1}{10}$ nonzero coefficients. Implementation details are available in~\ref{apx:implementation_details}. We chose the 3 datasets described in Table~\ref{tab:datasets}
\begin{table}
\caption{Description of datasets.} \label{tab:datasets}
\centering
\begin{tabular}{lrrrrr}
\toprule
{\bfseries\sffamily Dataset} & \multicolumn{1}{c}{$n$} & \multicolumn{1}{c}{$p$} & {\tablefont{density}} & \multicolumn{1}{c}{$L$} & $\Delta$\\
\midrule
{\bfseries\sffamily KDD 2010}~\citep{yu2010feature} & \hfill 19,264,097 & \hfill 1,163,024 & \hfill $10^{-6}$ & \hfill 28.12 & 0.15\\
{\bfseries\sffamily KDD 2012}~\citep{juan2016field} & \hfill 149,639,105 & \hfill 54,686,452 & \hfill $2 \times 10^{-7}$ & \hfill $1.25$ & 0.85\\
{\bfseries\sffamily Criteo}~\citep{juan2016field} & \hfill 45,840,617 & \hfill 1,000,000 & \hfill $4 \times 10^{-5}$ & \hfill $1.25$ & 0.89\\
\bottomrule
\end{tabular}
\end{table}

\begin{figure*}
\makebox[\textwidth][c]{
\includegraphics[width=1.0\linewidth]{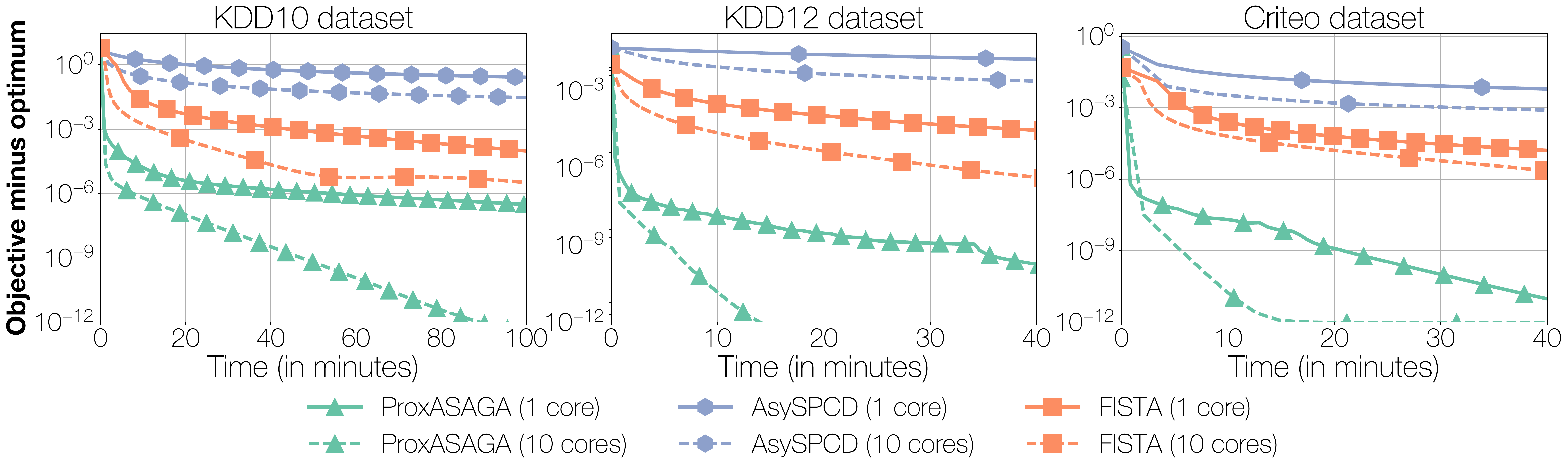}
}
\makebox[\textwidth][c]{
\includegraphics[width=1.0\linewidth]{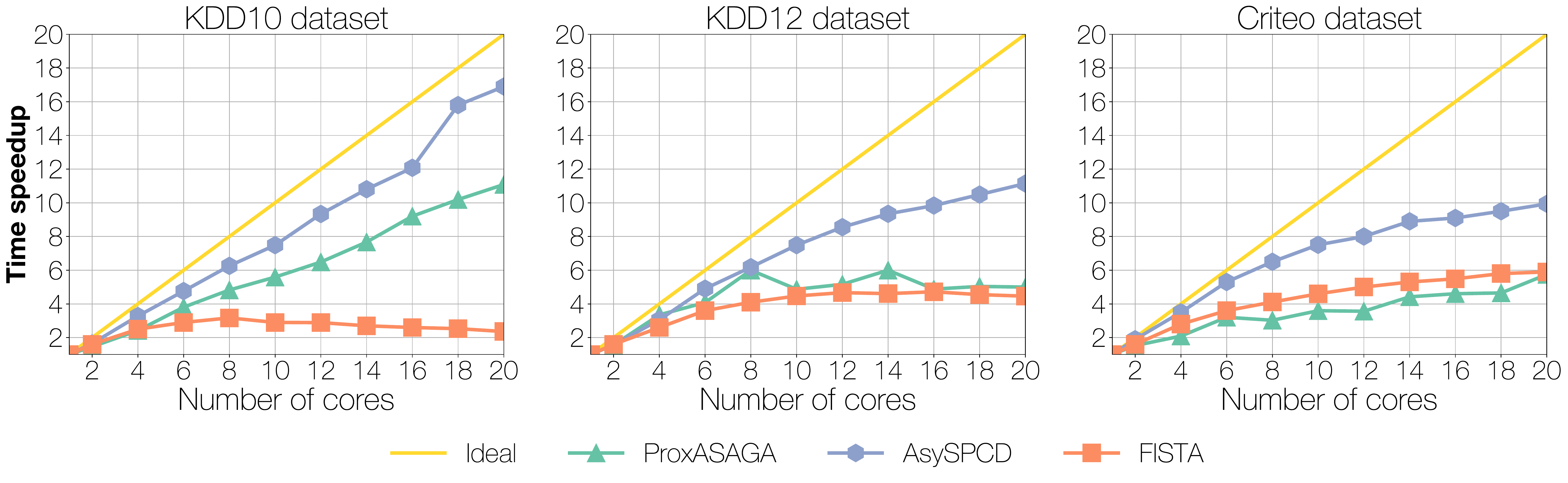}
}
\caption{{\bfseries Convergence for asynchronous stochastic methods for $\ell_1 + \ell_2$-regularized logistic regression}. {Top}: Suboptimality as a function of time for different asynchronous methods using 1 and 10 cores.
{Bottom}: Running time speedup as function of the number of cores. \PASAGA\ achieves significant speedups over its sequential version while being orders of magnitude faster than competing methods.
\AsySPCD\ achieves the highest speedups but it also the slowest overall method.
}\label{fig:suboptimality}
\end{figure*}

\paragraph{Results.} We compare three parallel asynchronous methods on the aforementioned datasets: \PASAGA\ (this work),\footnote{A reference C++/Python implementation of is available at \url{https://github.com/fabianp/ProxASAGA}} \AsySPCD, the asynchronous proximal coordinate descent method of~\citet{liu2015asynchronous2} and the (synchronous) \textsc{Fista} algorithm~\citep{beck2009gradient}, in which the gradient computation is parallelized by splitting the dataset into equal batches.
We aim to benchmark these methods in the most realistic scenario possible; to this end we use the following step size: $1 / 2L$ for \PASAGA, $1/L_c$ for \AsySPCD, where $L_c$ is the coordinate-wise Lipschitz constant of the gradient, while \textsc{Fista} uses backtracking line-search.
The results can be seen in Figure~\ref{fig:suboptimality} (top) with both one (thus sequential) and ten processors.
Two main observations can be made from this figure. First, \PASAGA\ is significantly faster on these problems. Second, its asynchronous version offers a significant speedup over its sequential counterpart.

In Figure~\ref{fig:suboptimality} (bottom) we present speedup with respect to the number of cores, where speedup is computed as the time to achieve a suboptimality of $10^{-10}$ with one core divided by the time to achieve the same suboptimality using several cores.
While our \emph{theoretical speedups} (with respect to the number of iterations) are almost linear as our theory predicts (see~\ref{apx:experiments}),  we observe a different story for our \emph{running time} speedups.
This can be attributed to memory access overhead, which our model does not take into account.
As predicted by our theoretical results, we observe a high correlation between the $\Delta$ dataset sparsity measure and the empirical speedup: KDD 2010 ($\Delta=0.15$) achieves a 11x speedup, while in Criteo ($\Delta=0.89$) the speedup is never above 6x.

Note that although competitor methods exhibit similar or sometimes better speedups, they remain orders of magnitude slower than \PASAGA\ in running time for large sparse problems. In fact, our method is between 5x and 80x times faster (in time to reach $10^{-10}$ suboptimality) than \textsc{Fista} and between 13x and 290x times faster  than \AsySPCD\ (see~\ref{apx:timing_benchmarks}).

\section{Conclusion and future work}
In this work, we have described \PASAGA, an asynchronous variance reduced algorithm with support for composite objective functions.
This method builds upon a novel sparse variant of the (proximal) \SAGA\ algorithm that takes advantage of sparsity in the individual gradients.
We have proven that this algorithm is linearly convergent under a condition on the step size and that it is linearly faster than its sequential counterpart given a bound on the delay.
Empirical benchmarks show that \PASAGA\ is orders of magnitude faster than existing state-of-the-art methods. %

This work can be extended in several ways.
First, we have focused on the \SAGA\ method as the basic iteration loop, but this approach can likely be extended to other proximal incremental schemes such as \SGD\ or \ProxSVRG.
Second, as mentioned in \S\ref{ssec:analysis}, it is an open question whether it is possible to obtain convergence guarantees without any sparsity assumption, as was done for \ASAGA.

\section*{Acknowledgements}

The authors would like to thank our colleagues Damien Garreau, Robert Gower, Thomas Kerdreux, Geoffrey Negiar and Konstantin Mishchenko for their feedback on this manuscript, and Jean-Baptiste Alayrac for support managing the computational resources.

This work was partially supported by a Google Research Award.
FP acknowledges support from the chaire \emph{\'Economie des nouvelles donn\'ees} with the \emph{data science} joint research initiative with the \emph{ fonds AXA pour la recherche}.

\bibliography{index}{}
\bibliographystyle{icml2017}

\clearpage
\appendix
\gdef\thesection{Appendix \Alph{section}}

\section*{\centering\LARGE Breaking the Nonsmooth Barrier: A Scalable Parallel Method for Composite Optimization\\
\hfill\\
\centering{Supplementary material}}

\vspace{2em}

\paragraph{Notations.} Throughout the supplementary material we use the following extra notation. We denote by $\langle \cdot, \cdot \rangle_{(i)}$ (resp. $\|\cdot\|_{(i)}$) the scalar product (resp. norm) restricted to blocks in $T_i$, i.e., $\langle \xx, \yy \rangle_{(i)} := \sum_{B \in T_i}\langle \llbracket\xx\rrbracket_B, \llbracket\yy\rrbracket_B \rangle$ and $\|\xx\|_{(i)} := \sqrt{\langle \xx, \xx\rangle_{(i)}}$.
We will also use the following definitions: $\varphi := \sum_{B \in \mathcal{B}} d_B h_B(\xx)$ and $\DD$ is the diagonal matrix defined block-wise as $\llbracket\DD\rrbracket_{B,B} = d_B \boldsymbol I_{|B|}$.

The {\bfseries Bregman divergence} associated with a convex function $f$ for points $\xx, \yy$ in its domain is defined as:
\begin{equation}
  B_f(\xx, \yy) := f(\xx) - f(\yy) - \langle \nabla f(\yy), \xx - \yy \rangle \,.
\end{equation}
  Note that this is always positive due to the convexity of $f$.

\section{Basic properties}

\begin{lemma}\label{lemma:neseterov_inequality}For any $\mu$-strongly convex function $f$ we have the following inequality:
\begin{equation}
\langle \nabla f(\yy) - \nabla f(\xx), \yy - \xx\rangle \geq \frac{\mu }{2}\|\yy - \xx\|^2  + B_f(\xx, \yy) \, .
\end{equation}

\end{lemma}
\begin{proof} By strong convexity, $f$ verifies the inequality:
  \begin{equation}\label{eq:strongconvex}
	f(\yy) \leq f(\xx) + \langle \nabla f(\yy), \yy - \xx \rangle - \frac{\mu}{2}\|\yy - \xx\|^2 \, ,
  \end{equation}

  for any $\xx, \yy$ in the domain (see e.g.~\citep{nesterov2004introductory}). We then have the equivalences:
\begin{align}
  &\qquad     f(\xx) \leq f(\yy) + \langle \nabla f(\xx), \xx - \yy \rangle - \frac{\mu}{2}\|\xx - \yy\|^2 \nonumber \\
  &\iff \frac{\mu}{2}\|\xx - \yy\|^2 + f(\xx) - f(\yy) \leq \langle \nabla f(\xx), \xx - \yy \rangle \nonumber \\
  &\iff \frac{\mu}{2}\|\xx - \yy\|^2 + \underbrace{f(\xx) - f(\yy) - \langle \nabla f(\yy), \xx - \yy\rangle}_{B_f(\xx, \yy)} \leq \langle \nabla f(\xx) - \nabla f(\yy), \xx - \yy \rangle \,,
\end{align}
  where in the last line we have subtracted $\langle \nabla f(\yy), \xx - \yy\rangle$ from both sides of the inequality.
\end{proof}

\hfill

\begin{lemma}\label{lemma:l_smooth_ineq}
  Let the $f_i$ be $L$-smooth and convex functions.
  Then it is verified that:
\begin{equation}
\frac{1}{n}\sum_{i=1}^n\|\nabla f_i(\xx) - \nabla f_i(\yy)\|^2 \leq 2 L B_f(\xx, \yy)\,.
\end{equation}
\end{lemma}

\begin{proof}
Since each $f_i$ is $L$-smooth, it is verified (see e.g.~\citet[Theorem 2.1.5]{nesterov2004introductory}) that
\begin{equation}
\|\nabla f_i(\xx) - \nabla f_i(\yy)\|^2 \leq 2 L \big(f_i(\xx) - f_i(\yy) - \langle \nabla f_i(\yy), \xx - \yy\rangle\big)\,.
\end{equation}
The result is obtained by averaging over $i$.
\end{proof}

\hfill

\begin{lemma}[Characterization of the proximal operator]\label{lemma:charac_prox}
  Let $h$ be convex lower semicontinuous. Then we have the following characterization of the proximal operator:
  \begin{equation}\label{eq:prox_char}
    \zz = \prox_{\gamma h}(\xx)  \iff  \frac{1}{\gamma}(\xx - \zz) \in \partial h(\zz)\,.
  \end{equation}
\end{lemma}
\begin{proof}
  This is a direct consequence of the first order optimality conditions on the definition of proximal operator, see e.g.~\citep{beck2009gradient, nesterov20013gradient}.
\end{proof}

\hfill

\begin{lemma}[Firm non-expansiveness]\label{lemmma:block_nonexpansive}
  Let $\xx, \tilde{\xx}$ be two arbitrary elements in the domain of $\varphi_i$ and $\zz, \tilde{\zz}$ be defined as $\zz := \eprox_{\varphi_i}(\xx)$, $\tilde{\zz} := \eprox_{\varphi_i}(\tilde{\xx})$. Then it is verified that:
\begin{equation}
\langle \zz - \tilde{\zz}, \xx - \tilde{\xx} \rangle_{(i)} \geq \|\zz - \tilde{\zz}\|_{(i)}^2\,.
\end{equation}
\end{lemma}
\begin{proof}
  By the block-separability of $\varphi_i$, the proximal operator is the concatenation of the proximal operators of the blocks. In other words, for any block $B \in T_i$ we have:
\begin{equation}
  \llbracket\zz\rrbracket_B = \prox_{\gamma \varphi_B}(\llbracket\xx\rrbracket_B)~,\quad \llbracket\tilde{\zz}\rrbracket_B = \prox_{\gamma \varphi_B}(\llbracket\tilde{\xx}\rrbracket_B)\,,
\end{equation}
where $\varphi_B$ is the restriction of $\varphi_i$ to $B$.
  By firm non-expansiveness of the proximal operator (see e.g.~\citet[Proposition 4.2]{bauschke2011convex}) we have that:
  $$
  \langle \llbracket\zz\rrbracket_B - \llbracket\tilde{\zz}\rrbracket_B, \llbracket\xx\rrbracket_B - \llbracket\tilde{\xx}\rrbracket_B\rangle \geq \|\llbracket\zz\rrbracket_B - \llbracket\tilde{\zz}\rrbracket_B\|^2 \,.
  $$
Summing over the blocks in $T_i$ yields the desired result.
\end{proof}

\clearpage
\section{Sparse Proximal SAGA}\label{apx:sparse_saga}

\vspace{0.5em}

This Appendix contains all proofs for Section~\ref{scs:sparse_prox_saga}. The main result of this section is Theorem~\ref{theorem:rates_sparse_saga}, whose proof is structured as follows:
\begin{itemize}
\item We start by proving four auxiliary results that will be used later on in the proofs of both  synchronous and asynchronous variants. 
The first is the unbiasedness of key quantities used in the algorithm. 
The second is a characterization of the solutions of \eqref{eq:opt} in terms of $f$ and $\varphi$ (defined below) in Lemma~\ref{lemma:fixed_point}.
The third is a key inequality in Lemma~\ref{lemma:gradient_mapping_1} that relates the gradient mapping to other terms that arise in the optimization.
The fourth is an upper bound on the variance terms of the gradient estimator, relating it to the Bregman divergence of $f$ and the past gradient estimator terms.
\item In Lemma~\ref{lemma:lyapunov_inequality}, we define an upper bound on the iterates $\|\xx_t - \xx^*\|^2$, called a Lyapunov function, and prove an inequality that relates this Lyapunov function value at the current iterate with its value at the previous iterate.
\item Finally, in the proof of Theorem~\ref{theorem:rates_sparse_saga} we use the previous inequality in terms of the Lyapunov function to prove a geometric convergence of the iterates.
\end{itemize}

%

We start by proving the following unbiasedness result, mentioned in \S\ref{scs:sparse_prox_saga}.
\begin{lemma}\label{lemma:unbiasedness}
Let $\DD_i$ and $\varphi_i$ be defined as in \S\ref{scs:sparse_prox_saga}. Then it is verified that 
$\Econd \DD_i = \boldsymbol{I}_p$ and $\Econd\, \varphi_i = h$.
\end{lemma}
\begin{proof}
Let $B \in \mathcal{B}$ an arbitrary block. We have the following sequence of equalities:
\begin{align}
\Econd \llbracket\DD_i\rrbracket_{B, B} &= \frac{1}{n}\sum_{i=1}^n \llbracket\DD_i\rrbracket_{B, B} = \frac{1}{n}\sum_{i=1}^n d_B \mathds{1}\{B \in T_i\}\boldsymbol{I}_{|B|}\\
&=  \frac{1}{n}\sum_{i=1}^n \frac{n}{n_B} \mathds{1}\{B \in T_i\}\boldsymbol{I}_{|B|}\\
&= \left(\frac{1}{n_B} \sum_{i=1}^n \mathds{1}\{B \in T_i\}\right) \boldsymbol{I}_{|B|}= \boldsymbol{I}_{|B|}~,
\end{align}
where the last equality comes from the definition of $n_B$. $\Econd \DD_i = \boldsymbol{I}_p$ then follows from the arbitrariness of $B$.

Similarly, for $\varphi_i$ we have:
\begin{align}
\Econd \varphi_i(\llbracket\xx\rrbracket_B) &= \frac{1}{n}\sum_{i=1}^n d_B \mathds{1}\{B \in T_i\}h_B(\llbracket\xx\rrbracket_B)\\
&=  \frac{1}{n}\sum_{i=1}^n \frac{n}{n_B} \mathds{1}\{B \in T_i\}h_B(\llbracket\xx\rrbracket_B)\\
&= \left(\frac{1}{n_B} \sum_{i=1}^n \mathds{1}\{B \in T_i\}\right) h_B(\llbracket\xx\rrbracket_B)= h_B(\llbracket\xx\rrbracket_B)~,
\end{align}
Finally, the result $\Econd\, \varphi_i = h$ comes from adding over all blocks.
\end{proof}

\begin{lemma}\label{lemma:fixed_point}
  $\xx^*$ is a solution to (OPT) if and only if the following condition is verified:
\begin{equation}
 \xx^* = \prox_{\gamma \varphi}\big(\xx^* - \gamma \boldsymbol \DD \nabla f(\xx^*)\big)\,.
\end{equation}
\end{lemma}
\begin{proof}By the first order optimality conditions, the solutions to \eqref{eq:opt} are characterized by the subdifferential inclusion ${-\nabla f(\xx^*) \in \partial h(\xx^*)}$.
  We can then write the following sequence of equivalences:
\begin{align}
    - \nabla f(\xx^*) \in \partial h(\xx^*) &\iff - \DD \nabla f(\xx^*) \in \DD \partial h(\xx^*)\nonumber\\
    &\qquad \text{ (multiplying by $\DD$, equivalence since diagonals are nonzero)}\nonumber\\
    &\iff - \DD \nabla f(\xx^*) \in \partial \varphi(\xx^*)\nonumber\\
    &\qquad \text{ (by definition of $\varphi$)}\nonumber\\
    &\iff \mfrac{1}{\gamma}(\xx^* - \gamma \DD \nabla f(\xx^*) - \xx^*)  \in \partial \varphi(\xx^*) \nonumber\\
    &\qquad \text{ (adding and subtracting $\xx^*$)}\nonumber\\
    &\iff \xx^* = \prox_{\gamma \varphi}(\xx^* - \gamma \DD \nabla f(\xx^*))\,.\\
    &\qquad \text{ (by Lemma \ref{lemma:charac_prox})} \nonumber
\end{align}
  Since all steps are equivalences, we have the desired result.
\end{proof}

\hfill

The following lemma will be key in the proof of convergence for both the sequential and the parallel versions of the algorithm.
With this result, we will be able to bound the product between the gradient mapping and the iterate suboptimality by:
\begin{itemize}
\item First, the negative norm of the gradient mapping, which will be key in the parallel setting to cancel out the terms arising from the asynchrony.
\item  Second, variance terms in $\|\vv_i - \DD_i \nabla f(\xx^*)\|^2$ that we will be able to bound by the Bregman divergence using Lemma~\ref{lemma:l_smooth_ineq}.
\item Third and last, a product with terms in $\langle \vv_i - \DD_i\nabla f(\xx^*) , \xx - \xx^* \rangle$, which taken in expectation gives $\langle \nabla f(\xx) - \nabla f(\xx^*) , \xx - \xx^* \rangle$ and will allow us to apply Lemma~\ref{lemma:neseterov_inequality} to obtain the contraction terms needed to obtain a geometric rate of convergence.
\end{itemize}

\hfill

\begin{lemma}[Gradient mapping inequality]\label{lemma:gradient_mapping_1}
  Let $\xx$ be an arbitrary vector, $\xx^*$ a solution to \eqref{eq:opt},
  $\vv_i$ as defined in \eqref{eq:SPS} and $\boldsymbol g = \boldsymbol g(\xx, \vv_i, i)$ the gradient mapping defined in~\eqref{eq:def_grad_mapping}.  Then the following inequality is verified for any $\beta > 0$:
  \begin{equation}
  \langle \boldsymbol g, \xx - \xx^*\rangle\geq -\frac{\gamma}{2}(\beta - 2)\| \boldsymbol g\|^2 -  \frac{\gamma}{2\beta}\| \vv_i -  \DD_i \nabla f(\xx^*) \|^2 + \langle  \vv_i - \DD_i \nabla f(\xx^*) , \xx - \xx^*\rangle \,.
  \end{equation}
\end{lemma}
\begin{proof}

By firm non-expansiveness of the proximal operator (Lemma~\ref{lemmma:block_nonexpansive}) applied to $\zz = \prox_{\gamma \varphi_i}(\xx - \gamma \vv_i)$ and $\tilde{\zz} = \prox_{\gamma \varphi_i}(\xx^* - \gamma \DD  \nabla f(\xx^*))$ we have:
\begin{equation}
  \|\zz - \tilde\zz\|_{(i)}^2 - \langle \zz - \tilde\zz, \xx - \gamma \vv_i - \xx^* + \gamma \DD \nabla f(\xx^*) \rangle_{(i)} \leq 0~.
\end{equation}
By the \eqref{eq:SPS} iteration we have $\xx^+ = \zz$ and by Lemma~\ref{lemma:charac_prox} we have that $\llbracket\zz\rrbracket_{T_i} = \llbracket\xx^*\rrbracket_{T_i}$, hence the above can be rewritten as
  \begin{equation}\label{eq:nonexpansive}
    \|\xx^+ - \xx^*\|_{(i)}^2 - \langle \xx^+ - \xx^*, \xx - \gamma \vv_i - \xx^* + \gamma \DD \nabla f(\xx^*) \rangle_{(i)} \leq 0~.
  \end{equation}
  We can now write the following sequence of inequalities
\begin{align}
    &\langle \gamma \boldsymbol g, \xx - \xx^*\rangle = \langle \xx - \xx^+, \xx - \xx^* \rangle_{(i)} \qquad \text{ (by definition and sparsity of $g$)}\nonumber\\
    &\qquad=  \langle \xx - \xx^+ + \xx^* - \xx^*, \xx - \xx^* \rangle_{(i)} \nonumber\\
    &\qquad= \|\xx - \xx^*\|_{(i)}^2 -  \langle \xx^+ - \xx^*, \xx - \xx^* \rangle_{(i)} \nonumber\\
    &\qquad\geq \|\xx - \xx^*\|_{(i)}^2 - \langle \xx^+ - \xx^*, 2 \xx - \gamma \vv_i- 2 \xx^* + \gamma \DD \nabla f(\xx^*) \rangle_{(i)} + \|\xx^+ - \xx^*\|^2_{(i)} \\
    &\qquad\qquad \text{ (adding Eq.~\eqref{eq:nonexpansive})} \nonumber\\
    &\qquad= \|\xx - \xx^+\|_{(i)}^2 + \langle \xx^+ - \xx^*, \gamma \vv_i -  \gamma \DD \nabla f(\xx^*) \rangle_{(i)} \qquad \text{ (completing the square)} \nonumber\\
    &\qquad= \|\xx - \xx^+\|_{(i)}^2 + \langle \xx - \xx^*, \gamma \vv_i -  \gamma \DD \nabla f(\xx^*) \rangle_{(i)}- \langle \xx - \xx^+, \gamma \vv_i-  \gamma \DD \nabla f(\xx^*) \rangle_{(i)}  \nonumber\\
    &\qquad\qquad \text{ (adding and substracting $\xx$)} \nonumber\\
    &\qquad\geq \Big(1 - \frac{\beta}{2}\Big)\|\xx - \xx^+\|_{(i)}^2 -  \frac{\gamma^2}{2\beta}\| \vv_i -  \DD \nabla f(\xx^*) \|_{(i)}^2 + \gamma\langle  \vv_i -  \DD \nabla f(\xx^*) , \xx - \xx^*\rangle_{(i)} \nonumber\\
    &\qquad\qquad \text{(Young's inequality ${2 \langle a, b \rangle \leq \frac{\|a\|^2}{\beta} + \beta \|b\|^2}$, valid for arbitrary $\beta > 0$)} \nonumber\\
    &\qquad\geq \Big(1 - \frac{\beta}{2}\Big)\|\xx - \xx^+\|_{(i)}^2 -  \frac{\gamma^2}{2\beta}\| \vv_i -  \DD_i \nabla f(\xx^*) \|^2 + \gamma\langle  \vv_i -  \DD_i \nabla f(\xx^*) , \xx - \xx^*\rangle \nonumber\\
    &\qquad\qquad \text{ (by definition of $\DD_i$ and using the fact that $\vv_i$ is $T_i$-sparse)}\nonumber\\
    &\qquad= \Big(1 - \frac{\beta}{2}\Big)\|\gamma \boldsymbol g\|^2 -  \frac{\gamma^2}{2\beta}\| \vv_i -  \DD_i \nabla f(\xx^*) \|^2 + \gamma\langle  \vv_i -  \DD \nabla f(\xx^*) , \xx - \xx^*\rangle \,,
\end{align}
  where in the last inequality we have used the fact that $\boldsymbol g$ is $T_i$-sparse.
Finally, dividing by $\gamma$ both sides yields the desired result.
\end{proof}

\begin{lemma}[Upper bound on the gradient estimator variance]\label{lma:variance}
For arbitrary vectors $\xx$, $({\boldsymbol\alpha}_i)_{i=0}^n$, and $\vv_i$ as defined in \eqref{eq:SPS} we have:
\begin{equation}\label{eq:variance_terms}
    \Econd \|\vv_i - \DD_i \nabla f(\xx^*)\|^2  \leq 4 L B_f(\xx, \xx^*) + 2\Econd\|{\boldsymbol\alpha}_i - \nabla f_i(\xx^*) \|^2 \, .
\end{equation}

\end{lemma}
\begin{proof}
We will now bound the variance terms. For this we have:
\begin{align}
&\Econd\| \vv_i -  \DD \nabla f(\xx^*) \|_{(i)}^2 = \Econd\|\nabla f_i(\xx) - \nabla f_i(\xx^*)  + \nabla f_i(\xx^*) - {\boldsymbol\alpha}_i + \DD_i \overline{{\boldsymbol\alpha}} - \DD \nabla f(\xx^*)\|_{(i)}^2 \nonumber\\
&\qquad\leq 2 \Econd \|\nabla f_i(\xx) - \nabla f_i(\xx^*)\|^2 + 2 \Econd \|  \nabla f_i(\xx^*) - {\boldsymbol\alpha}_i - (\DD\nabla f(\xx^*) - \DD \overline{{\boldsymbol\alpha}} ) \|_{(i)}^2 \nonumber\\
&\qquad\qquad \text{ (by inequality $\|a + b\|^2 \leq 2 \|a\|^2 + 2 \|b\|^2$)}\nonumber\\
&\qquad= 2 \Econd \|\nabla f_i(\xx) - \nabla f_i(\xx^*)\|^2+ 2\Econd\|\nabla f_i(\xx^*) - {\boldsymbol\alpha}_i \|^2 \nonumber\\
&\qquad\qquad- 4 \Econd \langle \nabla f_i(\xx^*) - {\boldsymbol\alpha}_i, \DD \nabla f(\xx^*) - \DD \overline{{\boldsymbol\alpha}}  \rangle_{(i)} + 2\Econd \|\DD \nabla f(\xx^*) - \DD \overline{{\boldsymbol\alpha}} \|_{(i)}^2\,.\\
&\qquad\qquad \text{ (developing the square)}\nonumber
\end{align}
We will now simplify the last two terms in the above expression. For the first of the two last terms we have:
\begin{align}
  &- 4\Econd \langle \nabla f_i(\xx^*) - {\boldsymbol\alpha}_i, \DD \nabla f(\xx^*) - \DD \overline{{\boldsymbol\alpha}}  \rangle_{(i)} =- 4 \Econd \langle \nabla f_i(\xx^*) - {\boldsymbol\alpha}_i, \DD \nabla f(\xx^*) - \DD \overline{{\boldsymbol\alpha}}  \rangle\\
  &\qquad \quad \text{ (support of first term)} \nonumber\\
  &\qquad= - 4 \langle \nabla f(\xx^*) - \overline{{\boldsymbol\alpha}}, \DD\nabla f(\xx^*) - \DD\overline{{\boldsymbol\alpha}}  \rangle \nonumber\\
  &\qquad= - 4 \|\nabla f(\xx^*) - \overline{{\boldsymbol\alpha}} \|_{\DD}^2 \,.
\end{align}
Similarly, for the last term we have:
\begin{align}
  2 \Econd \|\DD\nabla f(\xx^*) - \DD \overline{{\boldsymbol\alpha}} \|_{(i)}^2  &= 2\Econd \langle \DD_i \nabla f(\xx^*) - \DD_i \overline{{\boldsymbol\alpha}}, \DD \nabla f(\xx^*) - \DD \overline{{\boldsymbol\alpha}} \rangle \nonumber\\
  &= 2 \langle \nabla f(\xx^*) - \overline{{\boldsymbol\alpha}}, \DD \nabla f(\xx^*) - \DD \overline{{\boldsymbol\alpha}} \rangle \nonumber\\
  &\qquad \text{ (using Lemma~\ref{lemma:unbiasedness})}\nonumber\\
  &= 2\|\nabla f(\xx^*) - \overline{{\boldsymbol\alpha}} \|_{\DD}^2 \,.
\end{align}
and so the addition of these terms is negative and can be dropped. In all, for the variance terms we have
\begin{align}
    \Econd \|\vv_i - \DD \nabla f(\xx^*)\|_{(i)}^2 &\leq 2 \Econd \|\nabla f_i(\xx) - \nabla f_i(\xx^*)\|^2 + 2\Econd\|{\boldsymbol\alpha}_i - \nabla f_i(\xx^*) \|^2 \nonumber\\
    &\leq 4 L B_f(\xx, \xx^*) + 2\Econd\|{\boldsymbol\alpha}_i - \nabla f_i(\xx^*) \|^2 \, . \qquad \text{ (by Lemma~\ref{lemma:l_smooth_ineq}) }
\end{align}
\end{proof}

We now define an upper bound on the quantity that we would like to bound, often called a Lyapunov function, and establish a recursive inequality on this Lyapunov function.

\begin{lemma}[Lyapunov inequality]\label{lemma:lyapunov_inequality}
  Let $\mathcal{L}$ be the following $c$-parametrized function:
\begin{equation}\label{eq:lyapunov}
  \mathcal{L}(\xx,{\boldsymbol\alpha}) := \|\xx - \xx^*\|^2 + \frac{c}{n} \sum_{i=1}^n \|{\boldsymbol\alpha}_i - \nabla f_i(\xx^*)\|^2 \, .
\end{equation}
  Let $\xx^+$ and ${\boldsymbol\alpha}^+$ be obtained from the Sparse Proximal \SAGA\ updates~\eqref{eq:SPS}.
  Then we have:
\begin{align}
  \Econd \mathcal{L}(\xx^+,{\boldsymbol\alpha}^+) - \mathcal{L}(\xx,{\boldsymbol\alpha}) \leq  &- \gamma \mu \|\xx - \xx^*\|^2 +  \left(4 L \gamma^2 - 2 \gamma + 2L\frac{c}{n}\right) B_f(\xx, \xx^*)
  \nonumber \\
  &+ \left(2 \gamma^2 - \frac{c}{n}\right) \Econd \|{\boldsymbol\alpha}_i - \nabla f_i(\xx)\|^2\,.
\end{align}
\end{lemma}
\begin{proof}
  For the first term of $\mathcal{L}$ we have:
  \begin{align}
    \|\xx^+ - \xx^*\|^2 &= \|\xx - \gamma \boldsymbol g - \xx^*\|^2 \qquad \text{ ($\boldsymbol g := \boldsymbol g(\xx, \vv_i, i)$)}   \nonumber\\
    &= \|\xx - \xx^*\|^2 - 2\gamma \langle \boldsymbol g, \xx - \xx^* \rangle + \|\gamma \boldsymbol g\|^2 \nonumber\\
    &\leq \|\xx - \xx^*\|^2 +  \gamma^2\| \vv_i -  \DD_i \nabla f(\xx^*) \|^2 - 2 \gamma\langle  \vv_i - \DD_i \nabla f(\xx^*) , \xx - \xx^*\rangle \nonumber\\
    &\qquad \text{ (by Lemma~\ref{lemma:gradient_mapping_1} with $\beta=1$)}\nonumber\\
  \end{align}
Since $\vv_i$ is an unbiased estimator of the gradient and  $\Econd \DD_i =  \boldsymbol I_p$, taking expectations we have:
\begin{equation}\label{eq:proof_lyapunov_1}
  \begin{aligned}
    \Econd \|\xx^+ - \xx^*\|^2 &\leq \|\xx - \xx^*\|^2 +  \gamma^2\Econd\| \vv_i -  \DD_i \nabla f(\xx^*) \|^2 - 2 \gamma\langle  \nabla f(\xx) - \nabla f(\xx^*) , \xx - \xx^*\rangle\\
    &\leq (1 - \gamma \mu)\|\xx - \xx^*\|^2 +  \gamma^2\Econd\| \vv_i -  \DD_i \nabla f(\xx^*) \|^2 - 2 \gamma B_f(\xx, \xx^*)\,.\\
    &\qquad \text{ (by Lemma~\ref{lemma:neseterov_inequality})}
  \end{aligned}
\end{equation}

By using the variance terms bound (Lemma~\ref{lma:variance}) in the previous equation we have:
\begin{equation}\label{eq:lemma_lyapunov_3}
  \begin{aligned}
    \Econd \|\xx^+ - \xx^*\|^2
    &\leq (1 - \gamma \mu)\|\xx - \xx^*\|^2 + (4 L \gamma^2 - 2 \gamma) B_f(\xx, \xx^*)\\
    &\qquad + 2 \gamma^2 \Econd \|{\boldsymbol\alpha}_i - \nabla f_i(\xx^*)\|^2\,.\\
  \end{aligned}
\end{equation}

\hfill

We will now bound the second term of the Lyapunov function. We have:
  \begin{align}\label{eq:lemma_lyapunov_4}
    \frac{1}{n}\sum_{i=1}^n \|{\boldsymbol\alpha}_i^+ - \nabla f_i(\xx^*)\|^2 &= \left( 1 - \frac{1}{n}\right) \Econd \|{\boldsymbol\alpha}_i - \nabla f_i(\xx^*)\|^2 + \frac{1}{n}\Econd\|\nabla f_i(\xx) - \nabla f_i(\xx^*)\|^2  \nonumber\\
    &\qquad \text{ (by definition of $\balpha^+$)}\\
    &\leq \left( 1 - \frac{1}{n}\right) \Econd \|{\boldsymbol\alpha}_i - \nabla f_i(\xx^*)\|^2 + \frac{2}{n}L B_f(\xx, \xx^*)\,.\quad \text{ (by Lemma~\ref{lemma:l_smooth_ineq})}
  \end{align}
Combining Eq.~\eqref{eq:lemma_lyapunov_3} and \eqref{eq:lemma_lyapunov_4} we have:
\begin{align}
\Econd \mathcal{L}(\xx^+, {\boldsymbol\alpha}^+) &\leq (1 - \gamma \mu)\|\xx - \xx^*\|^2 +  (4 L \gamma^2 - 2 \gamma) B_f(\xx, \xx^*) + 2 \gamma^2 \Econd \|{\boldsymbol\alpha}_i - \nabla f_i(\xx^*)\|^2 \nonumber\\
 &\qquad + c\left[\left( 1 - \frac{1}{n}\right) \Econd \|{\boldsymbol\alpha}_i - \nabla f_i(\xx^*)\|^2 + \frac{1}{n}2 L B_f(\xx, \xx^*)\right] \nonumber\\
 &=  (1 - \gamma \mu)\|\xx - \xx^*\|^2 +  \left(4 L \gamma^2 - 2 \gamma + 2L\frac{c}{n}\right) B_f(\xx, \xx^*) \nonumber\\
  & \qquad + \left(2 \gamma^2 - \frac{c}{n}\right) \Econd \|{\boldsymbol\alpha}_i - \nabla f_i(\xx^*)\|^2 + c \Econd \|{\boldsymbol\alpha}_i - \nabla f_i(\xx^*)\|^2 \nonumber\\
  &= \mathcal{L}(\xx,{\boldsymbol\alpha})  - \gamma \mu \|\xx - \xx^*\|^2 + \left(4 L \gamma^2 - 2 \gamma + 2L\frac{c}{n}\right) B_f(\xx, \xx^*) \nonumber \\
  & \qquad+ \left(2 \gamma^2 - \frac{c}{n}\right) \Econd \|{\boldsymbol\alpha}_i - \nabla f_i(\xx^*)\|^2\,.
\end{align}
Finally, subtracting $\mathcal{L}(\xx,{\boldsymbol\alpha})$ from both sides yields the desired result.
\end{proof}

\begin{customtheorem}{1}
Let $\gamma = \frac{a}{5L}$ for any $a\leq 1$ and $f$ be $\mu$-strongly convex. Then Sparse Proximal \SAGA\ converges geometrically in expectation with a rate factor of at least $\rho = \frac{1}{5} \min\{\frac{1}{n}, a\frac{1}{\kappa}\}$. That is, for $\xx_t$ obtained after $t$ updates and $\xx^*$ the solution to \eqref{eq:opt}, we have the bound:
  $$
  \EE \|\xx_t - \xx^*\|^2 \leq (1-\rho)^t C_0\,,\quad \text{ with } C_0 :=  \|\xx_0 - \xx^*\|^2 + \textstyle\frac{1}{5 L^2} \sum_{i=1}^n\|{\boldsymbol\alpha}_i^0 - \nabla f_i(\xx^*)\|^2 \quad.
  $$
\end{customtheorem}
\begin{proof}
  Let $\overline{H} := \frac{1}{n}\sum_i \|{\boldsymbol\alpha}_i - \nabla f_i(\xx^*)\|^2$.
  By the Lyapunov inequality from Lemma~\ref{lemma:lyapunov_inequality}, we have:
\begin{align}
  &\Econd\mathcal{L}_{t+1} - (1 - \rho)\mathcal{L}_t \leq \rho \mathcal{L}_t - \gamma \mu\|\xx_t - \xx^*\|^2 +  \left(4 L \gamma^2 - 2 \gamma + 2L\frac{c}{n}\right)B_f(\xx_t, \xx^*) + \left(2\gamma^2 - \frac{c}{n} \right)\overline{H} \nonumber\\
  &\qquad=\left(\rho - \gamma \mu\right)\|\xx_t - \xx^*\|^2 +  \left(4 L \gamma^2 - 2 \gamma + 2L\frac{c}{n}\right)B_f(\xx_t, \xx^*) + \left[2\gamma^2 + c \left(\rho - \frac{1}{n}\right) \right]\overline{H} \nonumber\\
  &\qquad\qquad \text{ (by definition of $\mathcal{L}_t$)} \nonumber\\
  &\qquad\leq\left(\rho - \gamma \mu\right)\|\xx_t - \xx^*\|^2 +  \left(4 L \gamma^2 - 2 \gamma + 2L\frac{c}{n}\right)B_f(\xx_t, \xx^*) + \left(2\gamma^2 - \frac{2 c}{3 n} \right)\overline{H} \nonumber\\
  &\qquad\qquad\text{ (choosing $\rho \leq \frac{1}{3 n}$)}\nonumber\\
  &\qquad=\left(\rho - \gamma \mu\right)\|\xx_t - \xx^*\|^2 +  \left(10 L \gamma^2 - 2 \gamma \right)B_f(\xx_t, \xx^*) \nonumber\\
  &\qquad\qquad\text{ (choosing $\frac{c}{n} = 3 \gamma^2$)} \nonumber\\
  &\qquad\leq\left(\rho - \frac{a \mu}{5L} \right)\|\xx_t - \xx^*\|^2 \qquad\text{ (for all $\gamma = \frac{a}{5L}, a \leq 1$)} \nonumber\\
  &\qquad\leq0 \,.\qquad\text{ (for $\rho \leq \frac{a}{5} \cdot\frac{\mu}{L}$)}
\end{align}
And so we have the bound:
\begin{equation}
\Econd\mathcal{L}_{t+1} \leq \left(1 - \min\Big\{\frac{1}{3n}, \frac{a}{5}\cdot \frac{1}{\kappa}\Big\}\right)\mathcal{L}_t \leq \left(1 - \frac{1}{5}\min\Big\{\frac{1}{n}, a\cdot \frac{1}{\kappa}\Big\}\right)\mathcal{L}_t \,,
\end{equation}
where in the last inequality we have used the trivial bound $\frac{1}{3n} \leq \frac{1}{5n}$ merely for clarity of exposition.
Chaining expectations from $t$ to $0$ we have:
\begin{align}
  \EE \mathcal{L}_{t+1} &\leq \left(1 - \frac{1}{5}\min\Big\{\frac{1}{n}, a\cdot \frac{1}{\kappa}\Big\}\right)^{t+1}\mathcal{L}_0 \nonumber\\
  &= \left(1 - \frac{1}{5}\min\Big\{\frac{1}{n}, a\cdot \frac{1}{\kappa}\Big\}\right)^{t+1} \left(\|\xx_0 - \xx^*\|^2 + \frac{3 a^2}{5^2 L^2} \sum_{i=1}^n\|{\boldsymbol\alpha}_i^0 - \nabla f_i(\xx^*)\|^2 \right) \nonumber\\
  &\leq \left(1 - \frac{1}{5}\min\Big\{\frac{1}{n}, a\cdot \frac{1}{\kappa}\Big\}\right)^{t+1} \left(\|\xx_0 - \xx^*\|^2 + \frac{1}{5 L^2} \sum_{i=1}^n\|{\boldsymbol\alpha}_i^0 - \nabla f_i(\xx^*)\|^2 \right) \\
  &\qquad \text{ (since $a \leq 1$ and $\nicefrac{3}{5} \leq 1$)} \,.\nonumber
\end{align}
The fact that $\mathcal{L}_t$ is a majorizer of $\|\xx_{t} - \xx^*\|^2$ completes the proof.
\end{proof}

\clearpage

\section{ProxASAGA}\label{apx:async}

\vspace{0.5em}

In this Appendix we provide the proofs for results from Section~\ref{sec:pasaga}, that is Theorem~\ref{thm:convergence}  (the convergence theorem for \PASAGA) and Corollary~\ref{thm:corollary} (its speedup result).
\paragraph{Notation.} Through this section, we use the following shorthand for the gradient mapping: $\boldsymbol g_t := \boldsymbol g(\hat{\xx}_t, \hat{\vv}^t_{i_t}, i_t)$.

\subsection{Proof outline.}\label{ssec:proof}

As in the smooth case ($h=0$), we start by using the definition of $\xx_{t+1}$  in Eq.~\eqref{eq:definition} to relate the distance to the optimum in terms of its previous iterates:
\begin{equation}
\begin{aligned}
\|\xx_{t+1} - \xx^*\|^2 = &\|\xx_{t} - \xx^*\|^2  + 2 \gamma \langle \hat \xx_t - \xx_t, \boldsymbol g_t \rangle +\gamma^2 \|\boldsymbol g_t\|^2 - 2 \gamma \langle \hat \xx_t - \xx^*, \boldsymbol g_t\rangle \,  .
\end{aligned}
\end{equation}
However, in this case $\boldsymbol g_t$ is not a gradient estimator but a gradient mapping, so we cannot continue as is customary -- by using the unbiasedness of the gradient in the $\langle\hat \xx_t - \xx^*, \boldsymbol g_t\rangle$ term together with the strong convexity of~$f$ (see~\citet[Section 3.5]{leblond2016Asaga}).

To circumvent this difficulty, we derive a tailored inequality for the gradient mapping (Lemma~\ref{lemma:gradient_mapping_1} in~\ref{apx:sparse_saga}), which in turn allows us to use the classical unbiasedness and strong convexity arguments to get the following inequality:
\begin{align}
a_{t+1}
\leq (1 - \frac{\gamma \mu}{2}) a_t + \gamma^2 &\EE \|\boldsymbol g_t\|^2 - 2 \gamma \EE B_f(\hat{\xx}_t, \xx^*) +
\underbrace{ \gamma \mu \EE \|\hat{\xx}_t - \xx\|^2 +2 \gamma \EE \langle \boldsymbol g_t,  \hat{\xx}_t - \xx_t \rangle}_{\text{additional asynchrony terms}} \\
&\underbrace{+ \gamma^2(\beta - 2) \EE \|\boldsymbol g_t\|^2 + \frac{\gamma^2}{\beta} \EE \|\hat \vv^t_{i_t} - \DD_{i_t} \nabla f(\xx^*)\|^2}_{\text{additional proximal and variance terms}} \,, \nonumber
\end{align}
where $a_t := \EE \|\xx_t - \xx^*\|^2$.
Note that since $f$ is strongly convex, $B_f(\hat{\xx}_t, \xx^*) \geq \frac{\mu}{2}\|\hat \xx_t - \xx^*\|^2 $.

In the smooth setting, one first expresses the additional asynchrony terms as linear combinations of past gradient variance terms $(\EE \|\boldsymbol g_u\|^2)_{0 \leq u \leq t}$.
Then one crucially uses the negative Bregman divergence term to control the variance terms. However, in our current setting, we cannot relate the norm of the gradient mapping $\EE \|\boldsymbol g_t\|^2$ to the Bregman divergence (from which $h$ is absent).
Instead, we use the negative term $\gamma^2(\beta - 1)\EE\|\boldsymbol g_t\|^2$ to control all the $(\EE \|\boldsymbol g_u\|^2)_{0 \leq u \leq t}$ terms that arise from asynchrony.

The rest of the proof consists in:

$i)$ expressing the additional asynchrony terms as linear combinations of $(\EE \|\boldsymbol g_u\|^2)_{0 \leq u \leq t}$, following~\citet[Lemma 1]{leblond2016Asaga};

$ii)$ expressing the last variance term, $\|\hat \vv^t_{i_t} - D_{i_t} \nabla f(\xx^*)\|^2$, as a linear combination of past Bregman divergences (Lemma~\ref{lma:variance} in~\ref{apx:sparse_saga} and Lemma~2 from \citet{leblond2016Asaga});

$iii)$ defining a Lyapunov function, $\mathcal{L}_t := {\sum_{u=0}^t (1 - \rho)^{t-u} a_u}$, and proving that it is bounded by a contraction given conditions on the maximum step size and delay.

\vspace{0.5em}

\subsection{Detailed proof}

\begin{customtheorem}{2}[Convergence guarantee and rate of \PASAGA]
	Suppose $\tau \leq \frac{1}{10\sqrt{\Delta}}$.
	For any step size $\gamma = \frac{a}{L}$ with $a \leq \frac{1}{36} \min\{1, \frac{6 \kappa}{\tau}\}$, the inconsistent read iterates of Algorithm~\ref{alg:theoretical} converge in expectation at a geometric rate factor of at least: $\rho(a) = \frac{1}{5} \min \big\{\frac{1}{n},  a \frac{1}{\kappa}\big\},$
	i.e. $\EE \|\hat \xx_t-\xx^*\|^2 \leq (1-\rho)^t \,  \tilde C_0$, where $\tilde C_0$ is a constant independent of $t$ ($\approx \frac{n \kappa}{a}C_0$ with $C_0$ as defined in Theorem~\ref{th1}).
\end{customtheorem}
\begin{proof}
  In order to get an \textbf{initial recursive inequality}, we first unroll the (virtual) update:
  \begin{align}
    \|\xx_{t+1} - \xx^*\|^2 &= \|\xx_t - \gamma \boldsymbol g_t - \xx^*\|^2 = \|\xx_t - \xx^*\|^2 + \|\gamma \boldsymbol g_t\|^2 - 2 \gamma \langle \boldsymbol g_t, \xx_t - \xx^*\rangle \nonumber \\
    &= \|\xx_t - \xx^*\|^2 + \|\gamma \boldsymbol g_t\|^2 - 2 \gamma \langle \boldsymbol g_t, \hat{\xx}_t - \xx^*\rangle + 2 \gamma \langle \boldsymbol g_t, \hat{\xx}_t - \xx_t\rangle \, ,
  \end{align}
  and then apply Lemma~\ref{lemma:gradient_mapping_1} with $\xx = \hat \xx_t$ and $\vv = \hat{\vv}_{i_t}^t$.
  Note that in this case we have $\boldsymbol g = \boldsymbol g_t$ and $\langle \cdot\rangle_{(i)} = \langle \cdot \rangle_{(i_t)}$.

  \begin{align}\label{eq:initrec}
	\|\xx_{t+1} - \xx^*\|^2
	&\leq \|\xx_t - \xx^*\|^2 + 2 \gamma \langle \boldsymbol g_t,  \hat{\xx}_t - \xx_t \rangle + \gamma^2\|\boldsymbol g_t\|^2 + \gamma^2(\beta - 2)\|\boldsymbol g_t\|^2\nonumber\\
	&\qquad  + \frac{\gamma^2}{\beta} \|\hat{\vv}_{i_t}^t - \DD \nabla f(\xx^*)\|_{({i_t})}^2
	- 2 \gamma \langle \hat{\vv}_{i_t}^t - \DD \nabla f(\xx^*), \hat{\xx}_t - \xx^* \rangle_{({i_t})} \nonumber\\
	&=\|\xx_t - \xx^*\|^2 + 2 \gamma \langle \boldsymbol g_t,  \hat{\xx}_t - \xx_t \rangle + \gamma^2(\beta - 1)\|\boldsymbol g_t\|^2\nonumber\\
	&\qquad  + \frac{\gamma^2}{\beta} \|\hat{\vv}_{i_t}^t - \DD_{i_t} \nabla f(\xx^*)\|^2
	- 2 \gamma \langle \hat{\vv}_{i_t}^t - \DD_{i_t} \nabla f(\xx^*), \hat{\xx}_t - \xx^* \rangle  . \\
	&
	\qquad \qquad \text{(as  ${[\hat{\vv}_{i_t}^t]_{T_{i_t}}} = \hat{\vv}_{i_t}^t$)} \nonumber
  \end{align}

  We now use the property that $i_t$ is independent of $\hat \xx_t$ (which we enforce by reading $\hat \xx_t$ before picking $i_t$, see Section~\ref{sec:pasaga}), together with the unbiasedness of the gradient update $\hat{\vv}_{i_t}^t$ ($\Econd \hat{\vv}_{i_t}^t = \nabla f(\hat{\xx}_t)$) and the definition of $\DD$ to simplify the following expression as follows:
  \begin{align}
    \Econd \langle \hat{\vv}_{i_t}^t - \DD_{i_t} \nabla f(\xx^*), \hat{\xx}_t - \xx^* \rangle &=  \langle \nabla f(\hat{\xx}_t) - \nabla f(\xx^*), \hat{\xx}_t - \xx^* \rangle \nonumber\\
    &\geq \frac{\mu }{2}\|\hat{\xx}_t - \xx^*\|^2  + B_f(\hat{\xx}_t, \xx^*) \, ,
  \end{align}
  where the last inequality comes from Lemma~\ref{lemma:neseterov_inequality}.
  Taking conditional expectations on~\eqref{eq:initrec} we get:
  \begin{align}
    \Econd\|\xx_{t+1} - \xx^*\|^2
    &\leq \|\xx_t - \xx^*\|^2 + 2 \gamma \Econd \langle \boldsymbol g_t,  \hat{\xx}_t - \xx_t \rangle + \gamma^2(\beta - 1) \Econd \|\boldsymbol g_t\|^2 \\
    &\qquad+ \frac{\gamma^2}{\beta} \Econd \|\hat{\vv}_{i_t}^t - \DD_{i_t} \nabla f(\xx^*)\|^2 - \gamma \mu\|\hat{\xx}_t - \xx^*\|^2 - 2 \gamma B_f(\hat{\xx}_t, \xx^*) \notag \\
    &\leq (1 - \frac{\gamma \mu}{2})\|\xx_t - \xx^*\|^2 + 2 \gamma \Econd \langle \boldsymbol g_t,  \hat{\xx}_t - \xx_t \rangle + \gamma^2(\beta - 1) \Econd\|\boldsymbol g_t\|^2\notag\\
    &\qquad + \frac{\gamma^2}{\beta} \Econd \|\hat{\vv}_{i_t}^t - \DD_{i_t} \nabla f(\xx^*)\|^2 \notag  + \gamma \mu\|\hat{\xx}_t - \xx_t\|^2  - 2 \gamma B_f(\hat{\xx}_t, \xx^*) \notag \\
    &\qquad \text{ (using $\|a+b\|^2 \leq 2 \|a\|^2 + 2\|b\|^2$ on $\|\xx_t - \hat{\xx}_t + \hat{\xx}_t - \xx^*\|^2$)} \notag \\
    &\leq (1 - \frac{\gamma \mu}{2})\|\xx_t - \xx^*\|^2 + \gamma^2(\beta - 1) \Econd\|\boldsymbol g_t\|^2 + \gamma \mu\|\hat{\xx}_t - \xx_t\|^2  + 2 \gamma \Econd \langle \boldsymbol g_t,  \hat{\xx}_t - \xx_t \rangle \notag \\
    &\qquad  - 2 \gamma B_f(\hat{\xx}_t, \xx^*) + \frac{4\gamma^2 L}{\beta}  B_f(\hat{\xx}_t, \xx^*) + \frac{2\gamma^2}{\beta} \Econd \|\hat{{\boldsymbol\alpha}}^t_{i_t} - \nabla f_{i_t}(\xx^*)\|^2 \,.\label{eq:BeginMaster} \\
    &\qquad \text{ (using Lemma~\ref{lma:variance} on the variance terms)} \notag
  \end{align}

  Since we also have:
  \begin{equation}
  \hat \xx_t - \xx_t = \gamma \sum_{u=(t - \tau)_+}^{t-1}\boldsymbol G_{u}^t \boldsymbol g(\hat \xx_{u}, \hat {\boldsymbol\alpha}^u, i_{u}),
  \end{equation}
  the effect of asynchrony for the perturbed iterate updates was already derived in a very similar setup in~\citet{leblond2016Asaga}.
  We re-use the following bounds from their Appendix C.4:\footnote{The appearance of the sparsity constant $\Delta$ is coming from the crucial property that $\Econd \|\xx\|_{(i)}^2 \leq \Delta \|\xx\|^2$ $\forall x \in \RR^p$ (see Eq.~(39) in~\citet{leblond2016Asaga}, where they use the notation $\|\cdot\|_i$ for our $\|\cdot\|_{(i)}$).}
  \begin{align}
    &\EE \|\hat{\xx}_t - \xx_t\|^2 \leq \gamma^2(1 + \sqrt{\Delta}\tau) \sum_{u=(t-\tau)_+}^{t-1} \EE \|\boldsymbol g_u\|^2 \,, & &\text{\citet[Eq.~(48)]{leblond2016Asaga}} \label{eq:XhatMinusXt}\\
    &\EE\langle \boldsymbol g_t, \hat{\xx}_t - \xx_t\rangle \leq \frac{\gamma \sqrt{\Delta}}{2} \sum_{u=(t-\tau)_+}^{t-1} \EE \|\boldsymbol g_u\|^2 + \frac{\gamma \sqrt{\Delta}\tau}{2}\EE \|\boldsymbol g_t\|^2 \, . & &\text{\citet[Eq.~(46)]{leblond2016Asaga}}   \label{eq:gtVsXhat}.
  \end{align}
  Because the updates on ${\boldsymbol\alpha}$ are the same for $\PASAGA$ as for $\ASAGA$, we can
  re-use the same argument arising in the proof of~\citet[Lemma 2]{leblond2016Asaga} to get the following bound on $\EE \|\hat{{\boldsymbol\alpha}}^t_{i_t} - \nabla f_{i_t}(\xx^*)\|^2$:
  \begin{equation} \label{eq:alphaitBound}
      \EE \|\hat{{\boldsymbol\alpha}}^t_{i_t} - \nabla f_{i_t}(\xx^*)\|^2 \leq \frac{2L}{n}\underbrace{\sum_{u=1}^{t - 1}(1 - \frac{1}{n})^{(t - 2 \tau - u  - 1)_+} \EE B_f(\hat{\xx}_u, \xx^*)}_{\text{\small Henceforth denoted $H_t$}} + 2 L (1 - \frac{1}{n})^{(t - \tau)_+}\tilde{e}_0 \,,
  \end{equation}
  where $\tilde e_0 := \frac{1}{2L} \EE\|{\boldsymbol\alpha}_i^0 - f'_i(\xx^*)\|^2$.
  This bound is obtained by analyzing which gradient could be the source of ${\boldsymbol\alpha}_{i_t}$ in the past (taking in consideration the inconsistent writes), and then applying Lemma~\ref{lemma:l_smooth_ineq} on the $\Econd \|\nabla f(\hat{\xx}_u) - \nabla f(\xx^*)\|^2$ terms, explaining the presence of $B_f(\hat{\xx}_u,\xx^*)$ terms.\footnote{Note that~\citet{leblond2016Asaga} analyzed the unconstrained scenario, and so
  $B_f(\hat{\xx}_u,\xx^*)$ is replaced by the simpler $f(\hat{x}_u) - f(\xx^*)$ in their bound.}
  The inequality~\eqref{eq:alphaitBound} corresponds to Eq.~(56) and~(57) in~\citet{leblond2016Asaga}.

  By taking the full expectation of~\eqref{eq:BeginMaster} and plugging the above inequalities back, we obtain an inequality similar to~\citet[Master inequality (28)]{leblond2016Asaga} which describes how the error terms $a_t := \EE \|\xx_t - \xx^*\|^2$ of the virtual iterates are related:
  \begin{equation}\label{eq:masterineq}
  \begin{aligned}
    a_{t+1} \leq
    &(1 - \frac{\gamma \mu}{2})a_t
       + \frac{4 \gamma^2 L}{\beta} (1 - \frac{1}{n})^{(t - \tau)_+}\tilde{e}_0 \\
    &+ \gamma^2\left[\beta - 1 + \sqrt{\Delta}\tau \right] \EE \|\boldsymbol g_t\|^2
       + \left[\gamma^2 \sqrt{\Delta} + \gamma^3 \mu (1 + \sqrt{\Delta}\tau) \right]\sum_{u=(t-\tau)_+}^t \EE \|\boldsymbol g_u\|^2\\
    &- 2 \gamma \EE B_f(\hat{\xx}_t, \xx^*) + \frac{4\gamma^2 L}{\beta} \EE B_f(\hat{\xx}_t, \xx^*)  + \frac{4 \gamma^2 L}{\beta n} H_t \, .
  \end{aligned}
  \end{equation}
  We now have a promising inequality with a contractive term and several quantities that we need to bound.
  In order to achieve our final result, we introduce the same Lyapunov function as in~\citet{leblond2016Asaga}:
  $$
  \begin{aligned}
    \mathcal{L}_t := \sum_{u=0}^t (1 - \rho)^{t-u} a_u \,,
  \end{aligned}
  $$
  where $\rho$ is a target rate factor for which we will provide a value later on.
  Proving that this Lyapunov function is bounded by a contraction will finish our proof.
  We have:
  \begin{align}
    \mathcal{L}_{t+1} = \sum_{u=0}^{t+1} (1 - \rho)^{t+1-u} a_u
    &= (1 - \rho)^{t+1}a_0
    + \sum_{u=1}^{t+1}(1 - \rho)^{t+1-u}a_u
    \nonumber \\
    &= (1 - \rho)^{t+1}a_0
    + \sum_{u=0}^t(1 - \rho)^{t-u}a_{u+1} \,  .
  \end{align}

  We now plug our new bound on $a_{t+1}$,~\eqref{eq:masterineq}:
  \begin{equation}
  \begin{aligned}
  \mathcal{L}_{t+1}
  \leq (1 - \rho)^{t+1} a_0
  + \sum_{u = 0}^t(1 - \rho)^{t-u} \Big[
  &(1 - \frac{\gamma\mu}{2}) a_u
  + \frac{4 \gamma^2 L}{\beta} (1 - \frac{1}{n})^{(u - \tau)_+}\tilde{e}_0 \\
  &+ \gamma^2 \big(\beta - 1 + \sqrt{\Delta}\tau \big) \EE \|\boldsymbol g_u\|^2 \\
  &+ \big(\gamma^2 \sqrt{\Delta} + \gamma^3 \mu (1 + \sqrt{\Delta}\tau) \big)\sum_{v=(u-\tau)_+}^u \EE \|\boldsymbol g_v\|^2\\
  &- 2 \gamma \EE B_f(\hat{\xx}_u, \xx^*) + \frac{4\gamma^2 L}{\beta} \EE B_f(\hat{\xx}_u, \xx^*)  + \frac{4 \gamma^2 L}{\beta n} H_u \Big] \,.
  \end{aligned}
  \end{equation}
  After regrouping similar terms, we get:
  \begin{equation}\label{eq:recLyap}
  \begin{aligned}
    \mathcal{L}_{t+1}
      \leq (1 - \rho)^{t+1} (a_0 + A \tilde{e}_0)
      + (1 - \frac{\gamma \mu}{2}) \mathcal{L}_t
      + \sum_{u=0}^t s_u^t \EE \|\boldsymbol g_u\|^2
      + \sum_{u=1}^t r_u^t  \EE B_f(\hat{\xx}_u, \xx^*) \, .
  \end{aligned}
  \end{equation}
  Now, provided that we can prove that under certain conditions the $s_u^t$ and $r_u^t$ terms are all negative (and that the $A$ term is not too big), we can drop them from the right-hand side of~\eqref{eq:recLyap} which will allow us to finish the proof.

  Let us compute these terms.
  Let $q := \frac{1-\nicefrac{1}{n}}{1-\rho}$ and we assume in the rest that $\rho < \nicefrac{1}{n}$.

  \textbf{Computing $A$.}
  We have:
  \begin{align}
    \frac{4 \gamma^2 L}{\beta} \sum_{u=0}^t (1 - \rho)^{t-u} (1 - \frac{1}{n})^{(u - \tau)_+}
    &\leq \frac{4 \gamma^2 L}{\beta} (1 - \rho)^t (1 - \rho)^{-\tau} (\tau +1+\frac{1}{1 - q}) \nonumber\\
    &\qquad \qquad \text{from~\citet[Eq (75)]{leblond2016Asaga}} \nonumber\\
    &= (1 - \rho)^{t+1} \underbrace{\frac{4 \gamma^2 L}{\beta} (1 - \rho)^{-\tau -1}(\tau +1+\frac{1}{1 - q})}_{:=A} \,.
  \end{align}

  \textbf{Computing $s_u^t$.}
  Since we have:
  \begin{equation}
    \sum_{u=0}^t (1 - \rho)^{t-u} \sum_{v = (u - \tau)_+}^{u-1}\EE \|\boldsymbol g_u\|^2
    \leq \tau (1 - \rho)^{-\tau} \sum_{u=0}^t (1 - \rho)^{t - u}  \EE\|\boldsymbol g_u\|^2 \,,
  \end{equation}
  we have for all $0 \leq u \leq t$:
  \begin{equation} \label{eq:sut}
  \begin{aligned}
  s_u^t \leq (1 - \rho)^{t - u} \Big[\gamma^2 \big(\beta -1 + \sqrt{\Delta} \tau) + \tau (1 - \rho)^{-\tau} \big(\gamma^2 \sqrt{\Delta} + \gamma^3 \mu (1 + \sqrt{\Delta}\tau) \big)\Big] \, .
  \end{aligned}
  \end{equation}

  \textbf{Computing $r_u^t$.}
  To analyze these quantities, we need to compute: $\sum_{u = 0}^t(1 - \rho)^{t-u} \sum_{v = 1}^{u-1} (1 - \frac{1}{n})^{(u - 2 \tau - v  - 1)_+}$.
  Fortunately, this is already done in~\citet[Eq (66)]{leblond2016Asaga}, and thus we know that for all $1\leq u \leq t$:
  \begin{equation} \label{eq:rut}
  \begin{aligned}
  r_u^t \leq (1 - \rho)^{t - u} \left[-2 \gamma + \frac{4\gamma^2 L}{\beta} + \frac{4L \gamma^2}{n \beta} (1 - \rho)^{-2 \tau -1} \Big( 2\tau+  \frac{1}{1 - q}\Big)\right]\,,
  \end{aligned}
  \end{equation}
  recalling that $q := \frac{1 - \nicefrac{1}{n}}{1 - \rho}$ and that we assumed $\rho < \frac{1}{n}$.

  We now need some assumptions to further analyze these quantities.
  We make simple choices for simplicity, though a tighter analysis is possible.
  To get manageable (and simple) constants, we follow~\citet[Eq.~(82) and (83)]{leblond2016Asaga} and assume:
  \begin{equation} \label{eq:initialRhoTauConditions}
  \begin{aligned}
  \rho \leq \frac{1}{4 n}  ;\quad \tau \leq \frac{n}{10} \,.
  \end{aligned}
  \end{equation}
  This tells us:
  \begin{equation*}
  \begin{aligned}
  \frac{1}{1 - q} &\leq \frac{4 n}{3}
  \\
  (1 - \rho)^{- k \tau - 1} &\leq \frac{4}{3} \qquad \qquad \text{for} \, \, 0 \leq k \leq 2 \,.
  \qquad \text{(using Bernouilli's inequality)}
  \end{aligned}
  \end{equation*}
  Additionally, we set $\beta = \frac{1}{2}$.
  Equation~\eqref{eq:sut} thus becomes:
  \begin{equation}
  s_u^t \leq \gamma^2 (1 - \rho)^{t - u} \left[-\frac{1}{2} + \sqrt{\Delta} \tau + \frac{4}{3} \big( \sqrt{\Delta} \tau + \gamma \mu \tau (1 + \sqrt{\Delta} \tau) \big)\right] \,.
  \end{equation}

  We see that for $s_u^t$ to be negative, we need $\tau = \mathcal{O}(\frac{1}{\sqrt{\Delta}})$.
  Let us assume that $\tau \leq \frac{1}{10 \sqrt{\Delta}}$.
  We then get:
  \begin{equation}
  \begin{aligned}
  s_u^t \leq \gamma^2 (1 - \rho)^{t - u} \left[ - \frac{1}{2} + \frac{1}{10} + \frac{4}{30} + \gamma \mu \tau \frac{4}{3} \frac{11}{10} \right]\,.
  \end{aligned}
  \end{equation}
  Thus, the condition under which all $s_u^t$ are negative boils down to:
  \begin{equation} \label{eq:stepTauCondition}
  \gamma \mu \tau \leq \frac{2}{11}\,.
  \end{equation}

  Now looking at the $r_u^t$ terms given our assumptions, the inequality~\eqref{eq:rut} becomes:
  \begin{align}
  r_u^t &\leq (1 - \rho)^{t - u} \left[ -2 \gamma + 8 \gamma^2 L + \frac{8 \gamma^2 L}{n} \frac{4}{3}\big(\frac{n}{5} + \frac{4n}{3}\big)\right] \nonumber\\
  &\leq (1 - \rho)^{t - u} \big(-2 \gamma + 36 \gamma^2 L\big)\,.
  \end{align}

  The condition for all $r_u^t$ to be negative then can be simplified down to:
  \begin{equation}
  \gamma \leq \frac{1}{18L}\,.
  \end{equation}

  We now have a promising inequality for proving that our Lyapunov function is bounded by a contraction.
  However we have defined $\mathcal{L}_t$ in terms of the virtual iterate $\xx_t$, which means that our result would only hold for a given $T$ fixed in advance, as is the case in~\citet{mania2015perturbed}.
  Fortunately, we can use the same trick as in~\citet[Eq.~(97)]{leblond2016Asaga}: we simply add $\gamma B_f(\hat \xx_t, \xx^*)$ to both sides in~\eqref{eq:recLyap}.
  $r_t^t$ is replaced by $r_t^t + \gamma$, which makes for a slightly worse bound on $\gamma$ to ensure linear convergence:
  \begin{equation}  \label{eq:stepCondition}
  \gamma \leq \frac{1}{36L} \, .
  \end{equation}
  For this small cost, we get a contraction bound on $B_f(\hat \xx_t, \xx^*)$, and thus by the strong convexity of $f$ (see~\eqref{eq:strongconvex}) we get a contraction bound for $\EE \|\hat \xx_t - \xx^*\|^2$.

  \textbf{Recap.}
  Let us use $\rho = \frac{1}{4 n}$ and $\gamma := \frac{a}{L}$. Then the conditions~\eqref{eq:stepTauCondition} and~\eqref{eq:stepCondition} on the step size~$\gamma$ reduce to:
  \begin{equation} \label{eq:FinalStepsizeCondition}
  a \leq \frac{1}{36} \min \{1, \frac{72}{11} \frac{\kappa}{\tau} \} .
  \end{equation}
  Moreover, the condition:
  \begin{equation} \label{eq:FinalTauCondition}
  \tau \leq \frac{1}{10 \sqrt{\Delta}}
  \end{equation}
  is sufficient to also ensure that~\eqref{eq:initialRhoTauConditions} is satisfied as $\Delta \in [\frac{1}{n}, 1]$, and thus $\frac{1}{\sqrt{\Delta}} \leq \sqrt{n} \leq n$.

  Thus under the conditions~\eqref{eq:FinalStepsizeCondition} and~\eqref{eq:FinalTauCondition},
  we have that all $s_u^t$ and $r_u^t$ terms are negative and we can rewrite the recurrent step of our Lyapunov function as:
  \begin{equation}\label{eq:finish}
  \mathcal{L}_{t+1} \leq \gamma \EE B_f(\hat{\xx}_t) + \mathcal{L}_{t+1}
  \leq (1 - \rho)^{t+1} (a_0 + A \tilde{e}_0)
  + (1 - \frac{\gamma \mu}{2}) \mathcal{L}_{t}\,.
  \end{equation}

  By unrolling the recursion~\eqref{eq:finish}, we can carefully combine the effect of the geometric term~$(1-\rho)$ with the one of~$(1-\frac{\gamma \mu}{2})$. This was already done in~\citet[Apx C.9, Eq. (101) to (103)]{leblond2016Asaga}, with a trick to handle various boundary cases, yielding the overall rate:
   \begin{equation}
	  \EE B_f(\hat \xx_t, \xx^*) \leq (1 - \rho^*)^{t+1} \hat{C}_0,
  \end{equation}
  where $\rho^* = \min\{\frac{1}{5n}, a \frac{2}{5\kappa}\}$ (that we simplified to $\rho^* = \frac{1}{5} \min\{\frac{1}{n}, a \frac{1}{\kappa}\}$ in the theorem statement). To get the final constant, we need to bound $A$.
  We have:
  \begin{align}
    A &= \frac{4 \gamma^2 L}{\beta} (1 - \rho)^{-\tau -1}(\tau +1+\frac{1}{1 - q})
    \nonumber\\
    &\leq 8 \gamma^2 L \frac{4}{3} (\frac{n}{10} + 1 + \frac{4n}{3})
	\nonumber    \\
    &\leq 26 \gamma^2 L n
    \nonumber\\
    &\leq \gamma n \, .
  \end{align}
  This is the same bound on $A$ that was used by~\citet{{leblond2016Asaga}} and so we obtain the same constant as their Eq.~(104):
  \begin{equation}
    \hat C_0 := \frac{21n}{\gamma}(\|\xx_0 - \xx^*\|^2 + \gamma \frac{n}{2L} \EE \|{\boldsymbol\alpha}_i^0 - \nabla f_i(\xx^*)\|^2) .
  \end{equation}
  Note that $\hat C_0 = \mathcal{O}(\frac{n}{\gamma} C_0)$ with $C_0$ defined as in Theorem~\ref{th1}.

  Now, using the strong convexity of~$f$ via~\eqref{eq:strongconvex}, we get:
  \begin{equation}
    \EE \|\hat \xx_t - \xx^*\|^2 \leq \frac{2}{\mu} \EE B_f(\hat \xx_t, \xx^*) \leq (1 - \rho^*)^{t+1} \tilde{C}_0,
  \end{equation}
  where $\tilde{C_0} = \mathcal{O}(\frac{n \kappa}{a} C_0)$.

  This finishes the proof for Theorem~\ref{thm:convergence}.
\end{proof}

\begin{customcorollary}{3}[Speedup]
	Suppose $\tau \leq \frac{1}{10 \sqrt{\Delta}}$.
	If $\kappa \geq n$, then using the step size $\gamma = \nicefrac{1}{36 L}$, \PASAGA\ converges geometrically with rate factor $\Omega(\frac{1}{\kappa})$.
	If $\kappa < n$, then using the step size $\gamma = \nicefrac{1}{36 n \mu}$, \PASAGA\ converges geometrically with rate factor $\Omega(\frac{1}{n})$.
	In both cases, the convergence rate is the same as Sparse Proximal \SAGA\ and \PASAGA\ is thus linearly faster than its sequential counterpart up to a constant factor. Note that in both cases \emph{the step size does not depend on~$\tau$}.

	Furthermore, if $\tau \leq 6 \kappa$, we can use a universal step size of $\Theta(\nicefrac{1}{L})$ to get a similar rate for \PASAGA\ than Sparse Proximal \SAGA, thus making it adaptive to local strong convexity since the knowledge of $\kappa$ is not required.
\end{customcorollary}
\begin{proof}
	If $\kappa \geq n$, the rate factor of Sparse Proximal \SAGA\ is $\nicefrac{1}{\kappa}$.
	To get the same rate factor, we need to choose $a = \Omega(1)$, which we can fortunately do since $\kappa \geq n \geq \sqrt{n} \geq 10 \frac{1}{10\sqrt{\Delta}} \geq 10 \tau$.

	If $\kappa < n$, then the rate factor of Sparse Proximal \SAGA\ is $\nicefrac{1}{n}$.
	Any choice of $a$ bigger than $\Omega(\nicefrac{\kappa}{n})$ gives us the same rate factor for \PASAGA.
	Since $\tau \leq \nicefrac{\sqrt{n}}{10}$ we can pick such an $a$ without violating the condition of Theorem~\ref{thm:convergence}.
\end{proof}

\clearpage

\section{Comparison with related work}\label{apx:related_work}

In this section, we relate our theoretical results and proof technique with the related literature.

\paragraph{Speedups.}
Our speedup regimes are comparable with the best ones obtained in the smooth case, including~\citet{hogwild2011,reddi2015variance}, even though unlike these papers, we support inconsistent reads and nonsmooth objective functions.
The one exception is~\citet{leblond2016Asaga}, where the authors prove that their algorithm, \ASAGA, can obtain a linear speedup even without sparsity in the well-conditioned regime.
In contrast, \PASAGA\ always requires some sparsity.
Whether this property for smooth objective functions could be extended to the composite case remains an open problem.

\paragraph{Coordinate Descent.}
We compare our approach for composite objective functions to its most natural competitor: \AsySPCD~\citep{liu2015asynchronous2}, an asynchronous stochastic coordinate descent algorithm.
While \AsySPCD\ also exhibits linear speedups, subject to a condition on $\tau$, one has to be especially careful when trying to compare these conditions.

First, while in theory the iterations of both algorithms have the same cost, in practice various tricks are introduced to save on computation, yielding different costs per updates.\footnote{For \PASAGA\ the relevant quantity becomes the average number of features per data point. For \AsySPCD\ it is rather the average number of data points per feature. In both cases the tricks involved are not covered by the theory.}
Second, the bound on $\tau$ for the coordinate descent algorithm depends on $p$, the dimensionality of the problem, whereas ours involves $n$, the number of data points.
Third, a more subtle issue is that $\tau$ is not affected by the same quantities for both algorithms.\footnote{To make sure $\tau$ is the same quantity for both algorithms, we have to assume that the iteration costs are homogeneous.}
See~\ref{apx:liu} for a more detailed explanation of the differences between the bounds.

In the best case scenario (where the components of the gradient are uncorrelated, a somewhat unrealistic setting), \AsySPCD\ can get a near-linear speedup for~$\tau$ as big as~$\sqrt[4]{p}$.
Our result states that $\tau = \mathcal{O}(\nicefrac{1}{\sqrt{\Delta}})$ is necessary for a linear speedup.
This means in case $\Delta \leq \nicefrac{1}{\sqrt{p}}$ our bound is better than the one obtained for \AsySPCD.
Recalling that $\nicefrac{1}{n} \leq \Delta \leq 1$, it appears that \PASAGA\ is favored when $n$ is bigger than $\sqrt{p}$ whereas \AsySPCD\ may have a better bound otherwise, though this comparison should be taken with a grain of salt given the assumptions we had to make to arrive at comparable quantities.

Furthermore, one has to note that while~\citet{liu2015asynchronous2} use the classical labeling scheme inherited from~\citet{hogwild2011}, they still assume in their proof that the $i_t$ are uniformly distributed and that their gradient estimators are conditionally unbiased -- though neither property is verified in the general asynchronous setting. Finally, we note that \AsySPCD\ (as well as its incremental variant Async-\textsc{ProxSvrcd}) assumes that the computation and assignment of the proximal operator is an atomic step, while we do not make such assumption.

\paragraph{SVRG.}
The Async-\ProxSVRG\ algorithm of~\citet{meng2017aaai} also exhibits theoretical linear speedups subject to the same condition as ours.
However, the analyzed algorithm uses dense updates and consistent read and writes.
Although they make the analysis easier, these two factors introduce costly bottlenecks and prevent linear speedups in running time.
Furthermore, here again the classical labeling scheme is used together with the unverified conditional unbiasedness condition.

\paragraph{Doubly stochastic algorithms.}
The Async-\textsc{ProxSvrcd} algorithm from~\citet{meng2017aaai,gu2016asynchronous} has a maximum allowable stepsize\footnote{To the best of our understanding, noting that extracting an interpretable bound from the given theoretical results was difficult.
	Furthermore, it appears that the proof technique may still have significant issues: for example, the ``fully lock-free'' assumption of~\citet{gu2016asynchronous} allows for overwrites, and is thus incompatible with their framework of analysis, in particular their Eq. (8).}
 that is in $\mathcal{O}(\nicefrac{1}{pL})$, whereas the maximum step size for \PASAGA\ is in $\Omega(\nicefrac{1}{L})$, so can be up to $p$ times bigger.
Consequently, \PASAGA\ enjoys much faster theoretical convergence rates.
Unfortunately, we could not find a condition for linear speedups to compare to.
We also note that their algorithm is not appropriate in a sparse features setting.
This is illustrated in an empirical comparison in~\ref{apx:experiments} where we see that their convergence in number of iterations is orders of magnitude slower than appropriate algorithms like \SAGA\ or \PASAGA.

\subsection{Comparison of bounds with~\citet{liu2015asynchronous2}}\label{apx:liu}
\paragraph{Iteration costs.}
For both \PASAGA\ and \AsySPCD, the average cost of an iteration is $\mathcal{O}(n \overline S)$ (where $\overline S$ is the average support size).
In the case of \PASAGA\ (see Algorithm~\ref{alg:theoretical}), at each iteration the most costly operation is the computation of $\overline {\boldsymbol\alpha}$, while in the general case we need to compute a full gradient for \AsySPCD.

In order to reduce these prohibitive computation costs, several tricks are introduced.
Although they lead to much improved empirical performance, it should be noted that in both cases these tricks are not covered by the theory.
In particular, the unbiasedness condition can be violated.

In the case of \PASAGA, we store the average gradient term $\overline {\boldsymbol\alpha}$ in shared memory.
The cost of each iteration then becomes the size of the extended support of the partial gradient selected at random at this iteration, hence it is in $\mathcal{O}(\Delta_l)$, where $\Delta_l := \max_{i=1..n} |T_i|$.

For \AsySPCD, following~\citet{peng2016arock} we can store intermediary quantities for specific losses (e.g. $\ell_1$-regularized logistic regression).
The cost of an iteration then becomes the number of data points whose extended support includes the coordinate selected at random at this iteration, hence it is in $\mathcal{O}(n \Delta)$.

The relative difference in update cost of both algorithms then depends heavily on the data matrix: if the partial gradients usually have a extended support but coordinates belong to few of them (this can be the case if $n \ll p$ for example), then the iterations of \AsySPCD\ can be cheaper than those of \PASAGA.
Conversely, if data points usually have small extended support but coordinates belong to many of them (which can happen when $p \ll n$ for example), then the updates of \PASAGA\ are the cheaper ones.

\paragraph{Dependency of $\tau$ on the data matrix.}
In the case of \PASAGA\ the sizes of the extended support of each data point are important -- they are directly linked to the cost of each iteration.
Identical iteration costs for each data point do not influence $\tau$, whereas heterogeneous costs may cause $\tau$ to increase substantially.
In contrast, in the case of \AsySPCD, the relevant parts of the data matrix are the number of data points each dimension touches -- for much the same reason.
In the bipartite graph between data points and dimensions, either the left or the right degrees matter for $\tau$, depending on which algorithm you choose.

In order to compare their respective bounds, we have to make the assumption that the iteration costs are homogeneous, which means that each data point has the same support size and each dimension is active in the same number of data points. This implies that $\tau$ is the same quantity for both algorithms.

\paragraph{Best case scenario bound for AsySPCD.}
The result obtained in~\citet{liu2015asynchronous2} states that if $\tau^2 \Lambda = \mathcal{O}(\sqrt{p})$, \AsySPCD\ can get a near-linear speedup (where $\Lambda$ is a measure of the interactions between the components of the gradient, with $1 \leq \Lambda \leq \sqrt{p}$).
In the best possible scenario where $\Lambda = 1$ (which means that the coordinates of the gradients are completely uncorrelated), $\tau$ can be as big as $\sqrt[4]{p}$.

\clearpage

\section{Implementation details}\label{apx:implementation_details}

\paragraph{Initialization.} In the Sparse Proximal \SAGA\ algorithm and its asynchronous variant, \PASAGA, the vector $\xx$ can be initialized arbitrarily. The memory terms $\balpha_i$ can be initialized to any vector that verifies $\text{supp}(\balpha_i) = \text{supp}(\nabla f_i)$. In practice we found that the initialization $\balpha_i = \boldsymbol 0$ is very fast to set up and often outperforms more costly initializations.

With this initialization, the gradient approximation before the first update of the memory terms becomes $\nabla f_i(\xx) + \DD_i \overline{\balpha}$. Since most of the values in $\balpha$ are zero, $\overline{\balpha}$ will tend to be small compared to $\nabla f_i(\xx)$, and so the gradient estimate is very close to the \SGD\ estimate $\nabla f_i(\xx)$. The \SGD\ approximation is known to have a very fast initial convergence (which, in light of Figure~\ref{fig:suboptimality}, our method inherits) and has even been used as a heuristic to use during the first epoch of variance reduced methods~\citep{schmidt2016minimizing}.

The initialization of coefficients $\xx_0$ was always set to zero.

\paragraph{Exact regularization.} Computing the gradient of a smooth regularization such as the squared $\ell_2$ penalty of Eq.~\eqref{eq:logistic_loss} is independent of $n$ and so we can use the exact regularizer in the update of the coefficients instead of storing it in $\balpha$, which would also destroy the compressed storage of the memory terms described below. In practice we use this ``exact regularization'', multiplied by $\DD_i$ to preserve the sparsity pattern. 

Assuming a squared $\ell_2$ regularization term of the form $\frac{\lambda}{2}$, the gradient estimate in~\eqref{eq:SPS} becomes (note the extra $\lambda \xx$)
\begin{equation}
\vv_i = \nabla f_i(\xx) - \balpha_i + \DD_i (\overline\balpha + \lambda \xx)~.
\end{equation}

\paragraph{Storage of memory terms.}
The storage requirements for this method is in the worst case a table of size $n \times p$.
However, as for \SAG\ and \SAGA, for linearly parametrized loss functions of the form $f_i(\xx) = \ell(\boldsymbol a_i^T \xx)$, where $\ell$ is some real-valued function and $(\boldsymbol a_i)_{i=1}^n$ are samples associated with the learning problem, 
this can be reduced to a table of size $n$~\citep[\S 4.1]{schmidt2016minimizing}. 
This includes popular linear models such as least squares or logistic regression with $\ell$ the squared or logistic function, respectively.

The reduce storage comes from the fact that in this case the partial gradients have the structure
\begin{equation}
\nabla f_i(\xx) = \boldsymbol a_i \underbrace{\ell'(\boldsymbol a_i^T \xx)}_{\text{ scalar }} \quad.
\end{equation}
Since $\boldsymbol a_i$ is independent of $\xx$, we only need to store the scalar $\ell'(\boldsymbol a_i^T \xx)$. This decomposition also explains why $\nabla f_i$ inherits the sparsity pattern of $\boldsymbol a_i$.

\paragraph{Atomic updates.} Most modern processors have support for atomic operations with minimal overhead. In our case, we implemented a double-precision atomic type using the {\verb!C++11!} atomic features (\verb!std::atomic<double>!). This type implements atomic operations through the compare and swap semantics.

Empirically, we have found it necessary to implement atomic operations at least in the vector $\balpha$ and $\overline\balpha$ to reach arbitrary precision. If non-atomic operations are used, the method converges only to a limited precision (around normalized function suboptimality of $10^{-3}$), which might be sufficient for some machine learning applications but which we found not satisfying from an optimization point of view.

\paragraph{AsySPCD.} Following~\citep{peng2016arock} we keep the vector $(\boldsymbol a_i^T \xx)_{i=1}^n$ in memory and update it at each iteration using atomic updates.

\paragraph{Hardware and software.}
All experiments were run on a Dell PowerEdge 920 machine with 4 Intel Xeon E7-4830v2
processors with 10 2.2GHz cores each and 384GB 1600 Mhz RAM.
The \PASAGA and \AsySPCD\ code was implemented on C++ and binded in Python. The \textsc{Fista} code is implemented in pure Python using NumPY and SciPy for matrix computations (in this case the bottleneck is in large sparse matrix-vector operations for which efficient BLAS routines were used). Our \PASAGA\ implementation can be downloaded from \url{http://github.com/fabianp/ProxASAGA}.

\clearpage

\section{Experiments}\label{apx:experiments}

All datasets used for the experiments were downloaded from the LibSVM dataset suite.\footnote{\url{https://www.csie.ntu.edu.tw/~cjlin/libsvmtools/datasets/}}

\subsection{\large Comparison of ProxASAGA with other sequential methods}

We provide a comparison between the Sparse Proximal \SAGA\ and related methods in the sequential case. We compare against two methods: the \textsc{Mrbcd} method of \citet{zhao2014accelerated} (which forms the basis of Async-\textsc{ProxSvrcd}) and the vanilla implementation of \SAGA~\citep{defazio2014saga}, which does not have the ability to perform sparse updates. 
We compare in terms of both passes through the data (epochs) and time.
We use the same step size for all methods $(1/ 3 L)$. Due to the slow convergence of some methods, we use a smaller dataset than the ones used in \S\ref{scs:experiments}. Dataset RCV1 has $n=697,641, d=47,236$ and  a density of $0.15$,
while Covtype is a dense dataset with $n=581,012, d=54$.

\begin{figure}[h]
  \centering\includegraphics[width=0.7\linewidth]{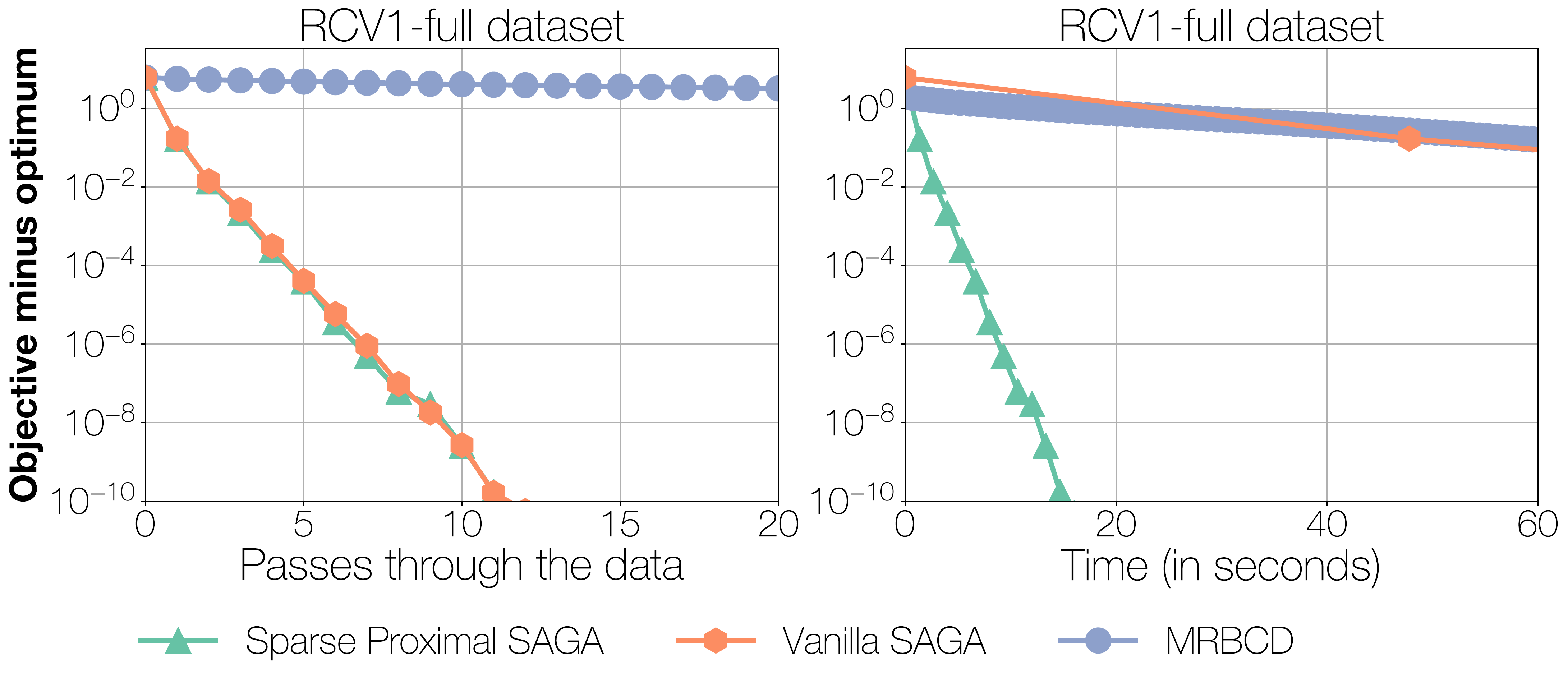}
  \centering\includegraphics[width=0.7\linewidth]{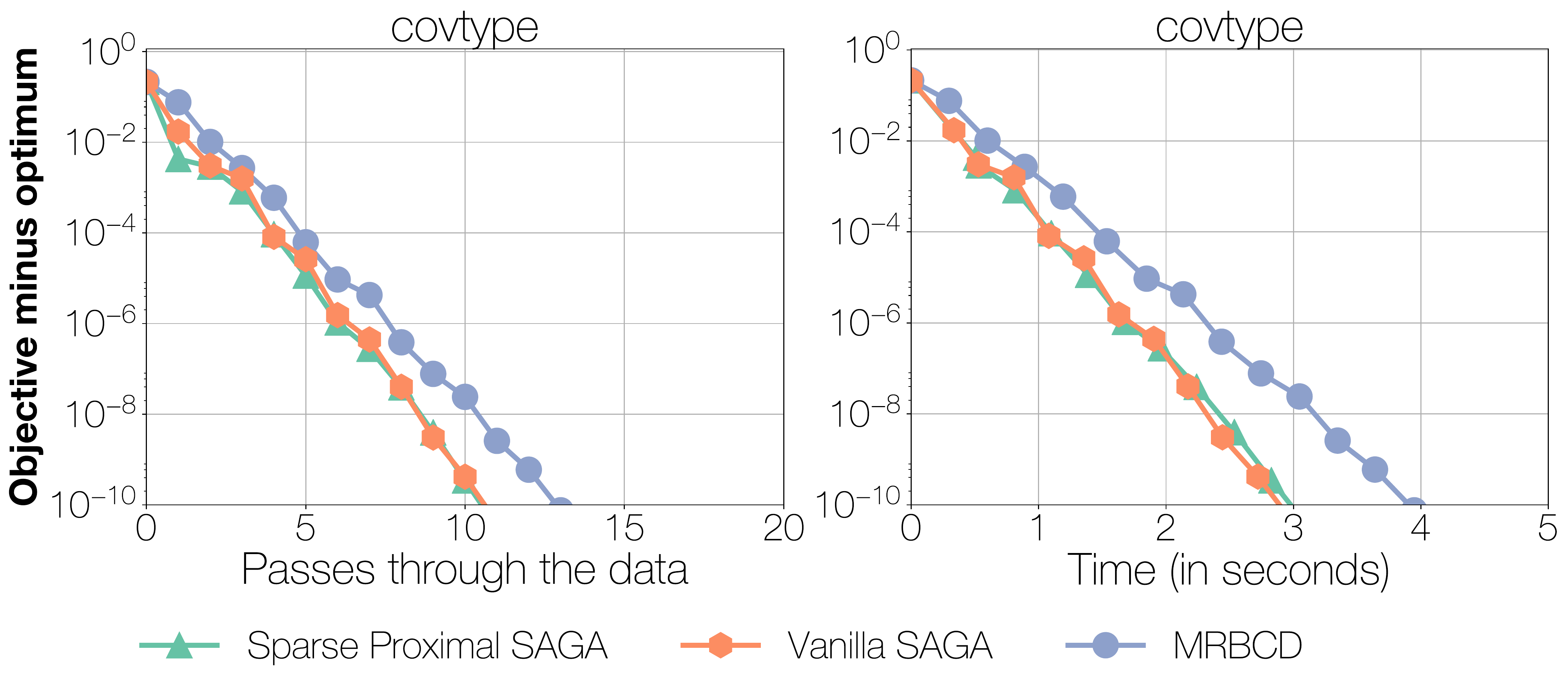}
  \caption{Suboptimality of different sequential algorithms. Each marker represents one pass through the dataset.}\label{fig:comparison_dense}
\end{figure}

We observe that for the convergence behavior in terms of number of passes, Sparse Proximal \SAGA\ performs as well as vanilla \SAGA, though the latter requires dense updates at every iteration (Fig.~\ref{fig:comparison_dense} top left). On the other hand, in terms of running time, our implementation of Sparse Proximal \SAGA\ is much more efficient than the other methods for sparse input (Fig.~\ref{fig:comparison_dense} top right). In the case of dense input (Fig.~\ref{fig:comparison_dense} bottom), the three methods perform similarly.

\vspace{-3mm}
\paragraph{A note on the performance of MRBCD.}
It may appear surprising that Sparse Proximal \SAGA\ outperforms \textsc{Mrbcd} so dramatically on sparse datasets.
However, one should note that \textsc{Mrbcd} is a doubly stochastic algorithm where both a random data point and a random coordinate are sampled for each iteration.
If the data matrix is very sparse, then the probability that the sampled coordinate is in the support of the sampled data point becomes very low.
This means that the gradient estimator term only contains the reference gradient term of \SVRG, which only changes once per epoch.
As a result, this estimator becomes very coarse and produces a slower empirical convergence.

This is reflected in the theoretical results given in~\citet{zhao2014accelerated}, where the epoch size needed to get linear convergence are $k$ times bigger than the ones required by plain \SVRG, where $k$ is the size of the set of blocks of coordinates.
\clearpage

\subsection{Theoretical speedups.}
In the experimental section, we have shown experimental speedup results where suboptimality was a function of the running time.
This measure encompasses both theoretical algorithmic optimization properties and hardware overheads (such as contention of shared memory) which are not taken into account in our analysis.

\begin{figure}[h]
	\includegraphics[width=\linewidth]{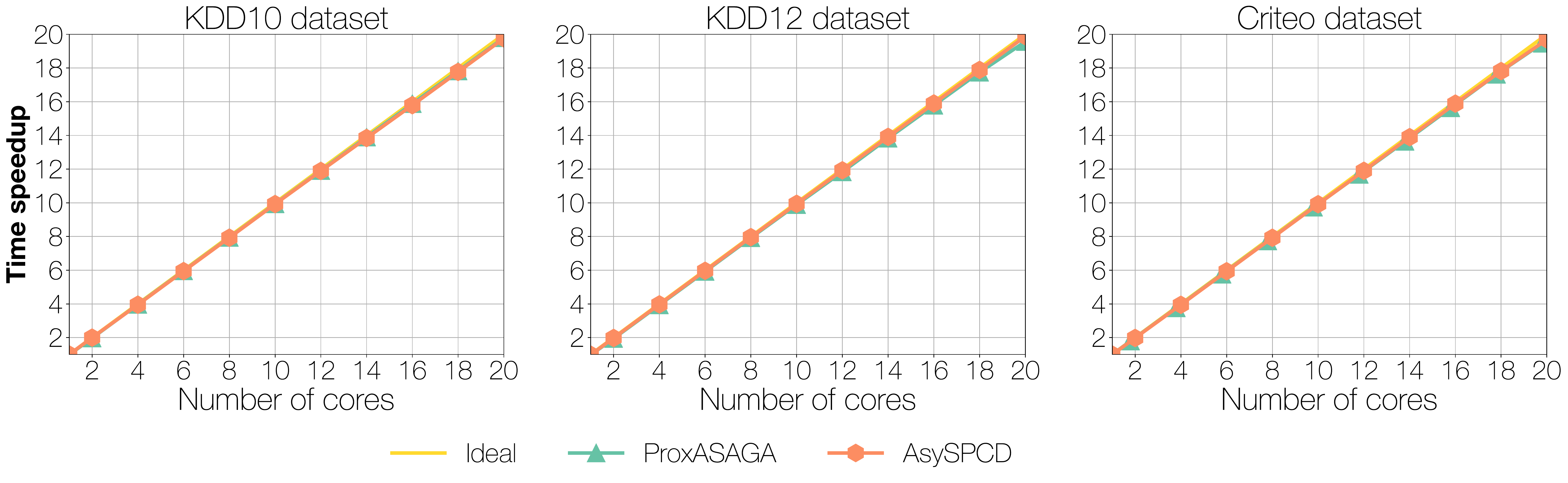}	
	\caption{{\bfseries Theoretical optimization speedups for $\ell_1 \!+\! \ell_2$-regularized logistic regression}. Speedup as measured by the number of iterations required to reach $10^{-5}$ suboptimality for \PASAGA\ and \AsySPCD. In \textsc{Fista} the iterates are the same with different cores and so matches the ``ideal'' speedup.
}\label{fig:ideal}
\end{figure}

In order to isolate these two effects, we now plot our speedup results in Figure~\ref{fig:ideal} where suboptimality is a function of the number of iterations; thus, we abstract away any potential hardware overhead.
	To do so, we implement a global counter which is sparsely updated (every $100$ iterations for example) in order not to modify the asynchrony of the system.
	This counter is used only for plotting purposes and is not needed otherwise.
Specifically, we define the theoretical speedup as:
\begin{equation*}
\text{theoretical speedup} := \text{(number of cores)} \, \frac{\text{number of iterations for sequential algorithm}}{\text{total number of iterations for parallel algorithm}} \, .
\end{equation*}

We see clearly that the theoretical speedups obtained by both \PASAGA and \AsySPCD\ are linear (i.e. ideal).
As we observe worse results in running time, this means that the hardware overheads of asynchronous methods are quite significant.

\subsection{Timing benchmarks}\label{apx:timing_benchmarks}
We now provide the time it takes for the different methods with 10 cores to reach a suboptimality of $10^{-10}$. All results are in hours.
\begin{table}[h]
\centering
\begin{tabular}{lrrr}
\toprule
{\bfseries\sffamily Dataset} & \multicolumn{1}{c}{\PASAGA} & \multicolumn{1}{c}{\AsySPCD} & \textsc{Fista} \\
\midrule
{\bfseries\sffamily KDD 2010} & \hfill 1.01 & \hfill 13.3 & \hfill 5.2 \\
{\bfseries\sffamily KDD 2012} & \hfill 0.09 & \hfill 26.6 & \hfill 8.3 \\
{\bfseries\sffamily Criteo} &   \hfill 0.14 & \hfill 33.3 & \hfill 6.6 \\
\bottomrule
\end{tabular}
\end{table}

\subsection{Hyperparameters}
The $\ell_1$-regularization parameter $\lambda_2$ was chosen as to give around 10\% of non-zero features. The exact chosen values are the following: $\lambda_2 = 10^{-11}$ for KDD 2010, $\lambda_2 = 10^{-16}$ for KDD 2012 and $\lambda_2 = 4 \times 10^{-12}$ for Criteo.

\clearpage

\end{document}